# SPECTRAL DISTRIBUTIONS AND ISOSPECTRAL SETS OF TRIDIAGONAL MATRICES

PETER C. GIBSON

ABSTRACT. We analyze the correspondence between finite sequences of finitely supported probability distributions and finite-dimensional, real, symmetric, tridiagonal matrices. In particular, we give an intrinsic description of the topology induced on sequences of distributions by the usual Euclidean structure on matrices. Our results provide an analytical tool with which to study ensembles of tridiagonal matrices, important in certain inverse problems and integrable systems. As an application, we prove that the Euler characteristic of any generic isospectral set of symmetric, tridiagonal matrices is a tangent number.

## Contents

## 1. INTRODUCTION

In this paper we analyze the correspondence between finite sequences of finitely supported probability distributions and finite-dimensional, real, symmetric, tridiagonal matrices. In particular, we give an intrinsic description of the topology induced

Research at MSRI was supported in part by NSF grant DMS-9810361. The author's research is supported also by NSERC Postdoctoral Fellowship 231108-2000.

on sequences of distributions by the usual Euclidean structure on matrices. Our results provide an analytical tool with which to study ensembles of tridiagonal matrices, which are important in a variety of contexts, including inverse problems [9] and integrable systems [1].

The paper is divided into two parts. In the first part, Section 2, we study the map $f$ that takes the spectral distribution of a finite Jacobi matrix to that matrix. The open $d$-simplex $S$ naturally represents probability distributions whose support is a fixed $(d+1)$-element set $\Lambda$ of reals, and it is well-known that $f$ maps $S$ homeomorphically onto the manifold $M$ of Jacobi matrices having spectrum $\Lambda$. However, viewed as a map on the closed simplex $\overline{S}$, $f$ has points of indeterminacy on the boundary $\partial S$. Our main result is the explicit construction of a blow-up $(B, \pi)$ of $\overline{S}$, centred at $\partial S$, that eliminates the points of indeterminacy of $f$, establishing a homeomorphism $h$ from $B$ onto the closure of $M$.

The original motivation was to better understand geometrically solutions to certain ill-posed inverse problems concerning Jacobi matrices, specifically [8, Theorem 6.4] and [9, Theorem 2.3]. The essential question, Problem 1 on p.7, is to describe intrinsically the topology induced on sequences of distributions by their correspondence with direct sums of Jacobi matrices. A complete answer to the question is implicit in the blow-up $(B, \pi)$ itself. As far as we know, these results are new; their application to inverse problems will appear in a future article.

The second part of the paper, Section 3, consists of an application of the first part to isospectral sets of real, symmetric, tridiagonal matrices. It turns out that $B$, although not convex, has a facial structure like that of a convex polytope. Based on this, we use the homeomorphism $h: B \to M$ to construct a family of polyhedral complexes $\overline{P}_d$ that encode the topology of the isospectral sets in question. This enables us to prove that the Euler characteristic of a generic isospectral set of symmetric, tridiagonal matrices is a tangent number. From the point of view of discrete geometry, the complexes $\overline{P}_d$ are remarkable in their own right. They have a simple description, yet exhibit exotic combinatorial properties.

The topology of isospectral sets of tridiagonal matrices has been considered independently, and much earlier, by Tomei, Nanda, Deift, Flaschka and others [10], [6], [7]. But whereas these authors make essential use of the Toda flow, our methods are more purely combinatorial. In place of a dynamical system, we exploit the abovementioned blow-up construction, which is completely explicit. Tomei also computed Euler characteristics, but does not seem to have noticed that hyperbolic tangent is a generating function for them. In any case, it is of interest to try to understand what connections exist between the two apparently different approaches.

*Acknowledgement.* Thanks both to Richard Cushman and Jürgen Bokowski for many helpful conversations at various points in the preparation of this article. Thanks also to Gunther Uhlmann and MSRI for the stimulating research environment of the MSRI semester in inverse problems, Fall 2001.

2. SEQUENCES OF DISTRIBUTIONS, DIRECT SUMS OF JACOBI MATRICES

Probability distributions proved useful in [9] for parametrizing solution sets of certain classes of inverse problems concerning finite Jacobi matrices, based on the well-known correspondence between distributions, orthogonal polynomials and Jacobi matrices—see [5]. In the present section, we study an extension of the standard correspondence, considering sequences of distributions and direct sums of Jacobi matrices in place of, respectively, distributions and Jacobi matrices. As will be explained in a future article, this extension gives a clearer geometric interpretation of the above-mentioned solution sets.

In Sections 2.1 and 2.2 we assemble the standard facts, without proof, and establish the basic notation to be used throughout Section 2. A precise formulation of the main problem is given at the end of Section 2.1. With the exception of the formula cited in Proposition 14 of Section 2.6, the material in Sections 2.3-2.8 is new and so we give detailed proofs of all results.

2.1. **Notation, basic facts, and definitions; Problem 1.** We are interested in the correspondence between probability distributions on the real line and Jacobi matrices, restricted to the finite-dimensional setting. With respect to this correspondence, sequences of orthogonal polynomials play a useful auxiliary role.

Let $\mathcal{J}$ denote the set of finite-dimensional Jacobi matrices, that is, matrices of the form

$$(1) \qquad J = \begin{pmatrix} a_1 & a_2 & 0 & \cdots & & 0 \\ a_2 & a_3 & a_4 & \ddots & & \vdots \\ 0 & a_4 & a_5 & \ddots & & 0 \\ \vdots & \ddots & \ddots & \ddots & & a_{2(d-1)} \\ 0 & \cdots & 0 & a_{2(d-1)} & a_{2d-1} \end{pmatrix},$$

for some $d \geq 1$, where each $a_n \in \mathbf{R}$, and each $a_n$ that has even index $n$ is strictly positive. We label the entries of $J$ as in (1) in order to suggest the map

$$(2) \qquad J \mapsto (a_1, a_2, \ldots, a_{2d-1}),$$

whereby $d \times d$ Jacobi matrices correspond to points in $\mathbf{R}^{2d-1}$.

Let $\mathcal{D}$ denote the set of finitely supported probability distributions on the real line. Elements of $\mathcal{D}$ are most simply represented by combinations of translates of the Dirac delta function, of the form

$$(3) \qquad d\alpha(\lambda) = \sum_{n=1}^{d} w_n \delta(\lambda - \lambda_n),$$

where each $\lambda_n \in \mathbf{R}$, each $w_n > 0$, and $\sum_{n=1}^{d} w_n = 1$. For technical reasons we shall adopt a slightly more cumbersome, equivalent representation, whereby the weights $w_n$ in (3) function as homogeneous, rather than absolute, coordinates. More precisely, let $\mathcal{D}'$ denote the set of distributions of the form (3), where each $\lambda_n, w_n \in \mathbf{R}$ and each $w_m w_n > 0$. For $d\alpha, d\beta \in \mathcal{D}'$, write
$$d\alpha \sim d\beta$$
if $d\alpha = r d\beta$ for some non-zero $r \in \mathbf{R}$. Then an element of $\mathcal{D}$ is represented by the equivalence class modulo $\sim$ of a distribution in $\mathcal{D}'$, consisting of its scalar multiples. In symbols, $\mathcal{D} = \mathcal{D}'/\sim$. Having established the precise meaning we will be less formal in practise; given $d\alpha \in \mathcal{D}'$, we will simply write $d\alpha \in \mathcal{D}$, letting the representative stand in for its class. The point of this notation is that we are free to rescale probability distributions without changing them.

Let $MOP$ denote the set of finite sequences of monic, orthogonal polynomials. These are sequences of the form
$$(p_0, p_1, \ldots, p_d),$$
where $d \geq 1$ and each $p_n(\lambda) \in \mathbf{R}[\lambda]$ is monic of degree $n$, and for which there exists a probability distribution $d\alpha$ on $\mathbf{R}$ such that
$$\langle p_n, p_n \rangle_{d\alpha} = \int_{-\infty}^{\infty} p_n^2 \, d\alpha \;>\; 0 \quad \text{for each } n,$$
$$\langle p_m, p_n \rangle_{d\alpha} = \int_{-\infty}^{\infty} p_m p_n \, d\alpha \;=\; 0 \quad \text{if } m \neq n.$$

Each $J \in \mathcal{J}$ corresponds to a unique $d\alpha \in \mathcal{D}$, its spectral distribution function, as follows. Orthogonally diagonalize $J$ as
$$J = O \Lambda O^T,$$
where $OO^T = I$ and $\Lambda = \mathrm{diag}(\lambda_1, \ldots, \lambda_d)$, and set

$$(4) \qquad d\alpha_J(\lambda) = \sum_{n=1}^{d} (O(1,n))^2 \, \delta(\lambda - \lambda_n),$$

where $O(i,j)$ denotes the $(i,j)$-entry of the matrix $O$. This is the (normalized representor of the) spectral distribution of $J$ in the sense that for any analytic function $f$,
$$f(J)(1,1) = \int_{-\infty}^{\infty} f \, d\alpha_J.$$

An equivalent characterization is the formula

$$(5) \qquad (zI - J)^{-1}(1,1) = \int_{-\infty}^{\infty} \frac{d\alpha_J(\lambda)}{z - \lambda}.$$

The map $J \mapsto d\alpha_J$ is a bijection from $\mathcal{J}$ onto $\mathcal{D}$. Our main focus in the present paper is its inverse
$$f : \mathcal{D} \to \mathcal{J}.$$
Here we define $f$ explicitly, along with an auxiliary map
$$\Psi : \mathcal{D} \to MOP.$$

**Definition 1.** *Let $d\alpha \in \mathcal{D}$ be given, and write $d = |\operatorname{supp} d\alpha|$. Define*

$$
(6) \qquad f(d\alpha) \;=\; \begin{pmatrix} a_1 & a_2 & 0 & \cdots & & 0 \\ a_2 & a_3 & a_4 & \ddots & & \vdots \\ 0 & a_4 & a_5 & \ddots & & 0 \\ \vdots & \ddots & \ddots & \ddots & & a_{2(d-1)} \\ 0 & \cdots & 0 & & a_{2(d-1)} & a_{2d-1} \end{pmatrix},
$$

$$
(7) \qquad \Psi(d\alpha) \;=\; (p_0, p_1, \ldots, p_d),
$$

*where the sequences $(a_1, \ldots, a_{2d-1})$ and $(p_0, \ldots, p_d)$ are constructed recursively as follows.*

<u>Initial step</u>  *For the purpose of notation set $a_0 = 0$. Set*

$$
p_0(\lambda) = 1, \quad a_1 = \frac{\langle 1, \lambda \rangle_{d\alpha}}{\langle 1, 1 \rangle_{d\alpha}}, \quad p_1(\lambda) = \lambda - a_1.
$$

<u>Continuing step</u>  *Given $a_{2n-2}, a_{2n-1}$ and $p_{n-1}, p_n$, for $n$ in the range $1 \le n \le d-1$, set*

$$
(8) \qquad a_{2n} \;=\; \sqrt{\frac{\langle p_n, p_n \rangle_{d\alpha}}{\langle p_{n-1}, p_{n-1} \rangle_{d\alpha}}},
$$

$$
(9) \qquad a_{2n+1} \;=\; \frac{\langle p_n, \lambda p_n \rangle_{d\alpha}}{\langle p_n, p_n \rangle_{d\alpha}},
$$

$$
(10) \qquad p_{n+1}(\lambda) \;=\; (\lambda - a_{2n+1}) p_n(\lambda) - a_{2n}^2 \, p_{n-1}(\lambda).
$$

The sequence of polynomials $\Psi(d\alpha)$ has a direct relationship to the matrix $J = f(d\alpha)$, namely, $p_n = \Psi_n(d\alpha)$ is the characteristic polynomial of the leading $n \times n$ submatrix of $J$; in particular, $p_d$ is the characteristic polynomial of $J$. (This fact will be needed in technical arguments later on.) Also, the truncated sequence $(p_0, \ldots, p_{d-1})$ is orthogonal with respect to the starting distribution $d\alpha$, but

$$
\int_{-\infty}^{\infty} p_d^2 \, d\alpha = 0
$$

is degenerate. Indeed $\operatorname{supp} d\alpha$ is the set of roots of $p_d$. If one extends $J$ to a $(d+1) \times (d+1)$ matrix $J'$ by choosing arbitrary entries $a_{2d}, a_{2d+1}$, then the full sequence $(p_0, \ldots, p_d)$ is orthogonal with respect to the spectral distribution $d\alpha_{J'}$ of $J'$. We index the components $f_n$ of $f$ according to the sequence $(a_1, \ldots, a_{2d-1})$, so that

$$
a_n = f_n(d\alpha).
$$

Let us consider a fixed $d = |\operatorname{supp} d\alpha|$, and write

$$
d\alpha(\lambda) = \sum_{n=1}^{d} w_n \delta(\lambda - \lambda_n).
$$

Observe that, according to Definition 1, for $1 \le n \le d-1$,

- $f_{2n}^2$ and $f_{2n+1}$ are rational functions in the variables $w_1, \ldots, w_d, \lambda_1, \ldots, \lambda_d$;
- $f_{2n}^2$ and $f_{2n+1}$ are homogeneous functions of degree 0 in the variables $w_1, \ldots, w_d$.

Now, $f$ itself is not rational, because of the square root (8). We have chosen the positive square root, and in some sense this choice is artificial. It is more natural to allow both square roots, and thus to make $f$ a multi-valued function (precisely, a $2^{d-1}$-valued function for fixed $d = |\operatorname{supp} d\alpha|$). This point of view is developed fully in Section 3; for now we leave $f$ as single-valued.

Having determined the nature of the functions $f_n$ in terms of the weights $w_n$, let us briefly consider the reverse question, restricted to to matrices $J$ having a prescribed set $\{\lambda_1, \ldots, \lambda_d\}$ of eigenvalues. To that end, let

$$\sum_{n=1}^{d} w_n \delta(\lambda - \lambda_n) = d\alpha_J \tag{11}$$

denote the normalized spectral distribution of $J$, as defined in (4). Then the normalized weights $w_n$ are polynomials in the entries $a_1, \ldots, a_{2d-1}$ of $J$, a fact easily inferred from (5), as follows. Let $J^{11}$ denote the submatrix of $J$ obtained by deleting the first row and column, and let

$$p_d(\lambda) = \prod_{n=1}^{d}(\lambda - \lambda_n)$$

be the characteristic polynomial of $J$. Then the left-hand side of (5) has the form

$$\det(zI - J^{11})/p_d(z),$$

while the right-hand side of (5), expressed in terms of (11), is the partial fraction decomposition

$$\sum_{n=1}^{d} \frac{w_n}{z - \lambda_n}.$$

Comparing these two expressions yields

$$w_n = \det(\lambda_n I - J^{11})/p_d'(\lambda_n),$$

which, for fixed $\lambda_1, \ldots, \lambda_d$, is a polynomial in the entries of $J$.

The set of $d \times d$ symmetric, tridiagonal matrices is a subspace of $\mathbf{R}^{d^2}$, isomorphic (indeed isometric) to $\mathbf{R}^{2d-1}$ by the map (2). With this isometry in mind, we will usually write $J \in \mathbf{R}^{2d-1}$ when $J$ is a $d \times d$ symmetric, tridiagonal matrix. Now, because $\mathcal{J} \cap \mathbf{R}^{2d-1}$ is not closed, it is useful to consider the extended class of matrices

$$\overline{\mathcal{J}} = \bigcup_{d \geq 1} \operatorname{cl}(\mathcal{J} \cap \mathbf{R}^{2d-1}).$$

$\overline{\mathcal{J}}$ consists of direct sums

$$J = J_1 \oplus \cdots \oplus J_r$$

of Jacobi matrices $J_1, \ldots, J_r \in \mathcal{J}$, where $r \geq 1$. By analogy we extend $\mathcal{D}$ to the class $\overline{\mathcal{D}}$ consisting of all sequences

$$d\alpha_1 \oplus \cdots \oplus d\alpha_r$$

of distributions $d\alpha_1, \ldots, d\alpha_r \in \mathcal{D}$, where $r \geq 1$. (We could alternatively write $(d\alpha_1, \ldots, d\alpha_r)$ for a sequence of distributions; the symbol $\oplus$ just serves to emphasize the connection with matrices.) The map

$$f : \mathcal{D} \to \mathcal{J}$$

extends to a map

$$\overline{f} : \overline{\mathcal{D}} \to \overline{\mathcal{J}}$$

in the obvious way, namely, by defining

$$\overline{f}(d\alpha_1 \oplus \cdots \oplus d\alpha_r) = f(d\alpha_1) \oplus \cdots \oplus f(d\alpha_r).$$

Note that, by its very definition, $\overline{f}$ is a bijection.

For each $d \geq 1$, $\overline{\mathcal{J}} \cap \mathbf{R}^{2d-1}$ inherits the Euclidean structure of $\mathbf{R}^{2d-1}$. But it is unclear a priori what it should mean for two points,

$$d\alpha_1 \oplus \cdots \oplus d\alpha_r, \quad d\beta_1 \oplus \cdots \oplus d\beta_s,$$

in $\overline{\mathcal{D}}$ to be "close". Thus, roughly speaking, a basic problem is to describe intrinsically the topology which $\overline{\mathcal{D}}$ acquires from $\overline{\mathcal{J}}$ via $\overline{f}$. We make this more specific as follows. Fix a $d$-element set $\Lambda \subset \mathbf{R}$. Denote by $\overline{\mathcal{J}}^\Lambda$ the set of $d \times d$ matrices $J \in \overline{\mathcal{J}}$ having spectrum $\Lambda$,

$$\overline{\mathcal{J}}^\Lambda = \left\{ J \in \overline{\mathcal{J}} \;\middle|\; J \in \mathbf{R}^{2d-1} \text{ and } \operatorname{spec} J = \Lambda \right\}.$$

Write $\overline{\mathcal{D}}^\Lambda$ to denote the corresponding set of sequences of distributions, that is, sequences

$$d\alpha_1 \oplus \cdots \oplus d\alpha_r \in \overline{\mathcal{D}}$$

such that $1 \leq r \leq d$,

(12) $$\sum_{n=1}^{r} |\operatorname{supp} d\alpha_n| = d \quad \text{and} \quad \bigcup_{n=1}^{r} \operatorname{supp} d\alpha_n = \Lambda.$$

With this notation, the restriction of $\overline{f}$,

$$\overline{f} : \overline{\mathcal{D}}^\Lambda \to \overline{\mathcal{J}}^\Lambda,$$

is a bijection.

> **Problem 1** Describe intrinsically the topology which $\overline{\mathcal{J}}^\Lambda$ induces on $\overline{\mathcal{D}}^\Lambda$ via $\overline{f}$. That is, characterize convergent sequences of points in $\overline{\mathcal{D}}^\Lambda$, and show how to compute limits of such sequences without evaluating $\overline{f}$.

We illustrate Problem 1 with a simple example. Let $d\alpha_0, d\alpha_1 \in \mathcal{D}$ denote the distributions

$$d\alpha_0(\lambda) = \delta(\lambda),$$
$$d\alpha_1(\lambda) = \delta(\lambda - 1) + \delta(\lambda - 2),$$

and consider the limit

$$L = \overline{f}^{-1}\left(\lim_{t\to 0^+} f(d\alpha_0 + t\,d\alpha_1)\right).$$

Given that the limit exists, what is $L$? A naive guess might be $L = d\alpha_0 \oplus d\alpha_1$, but this is wrong, as can be verified numerically by coding the recursive definition of $f$ (in, say, Matlab) and then computing $f(d\alpha_0 + td\alpha_t)$ for decreasing values of $t$. In fact, $L = d\alpha_0 \oplus d\alpha_1'$, where

$$\begin{aligned} d\alpha_1'(\lambda) &= \delta(\lambda - 1) + 4\delta(\lambda - 2) \\ &= \lambda^2 d\alpha_1(\lambda). \end{aligned}$$

In other words, for $\Lambda = \{0, 1, 2\}$, the curve $d\beta_t = d\alpha_0 + td\alpha_1$ (or equivalently, sequence, if we set $t = 1/n$, $n \in \mathbf{Z}^+$) converges in $\overline{\mathcal{D}}^\Lambda$, with respect to the topology induced by $\overline{f}$, to $d\alpha_0 \oplus \lambda^2 \, d\alpha_1$. And although $d\beta_t \in \mathcal{D}$ for every $t > 0$, its limiting value as $t \to 0^+$ is not in $\mathcal{D}$. Rather, $d\beta_t$ "splits" in the limit into a pair of smaller distributions. Problem 1 calls, in part, for a full description of this phenomenon.

2.2. **Preliminary observations, further notation.** Let $\Lambda \subset \mathbf{R}$ be a $d$-element set, for some $d \geq 1$, to be held fixed throughout the present section. As a first step toward solving Problem 1, we note that $\overline{\mathcal{D}}^\Lambda$ is a disjoint union of open discs of various dimensions. The essence of Problem 1 then consists in describing how these component discs fit together.

The notion of an ordered partition is useful for present considerations. Given a set $S$, a sequence

$$(S_1, \ldots, S_r)$$

of sets is an ordered partition of $S$ if: each $S_n$ is non-empty; $S_m \cap S_n = \emptyset$ whenever $m \neq n$; and $\bigcup_{n=1}^r S_n = S$. Denote the set of ordered partitions of $S$ by $\mathcal{P}_S$. The condition (12) for $d\alpha_1 \oplus \cdots \oplus d\alpha_r$ to belong to $\overline{\mathcal{D}}^\Lambda$ is precisely that

$$(\operatorname{supp} d\alpha_1, \ldots, \operatorname{supp} d\alpha_r)$$

be an ordered partition of $\Lambda$. And for each ordered partition $(S_1, \ldots, S_r) \in \mathcal{P}_\Lambda$, there corresponds a piece of $\overline{\mathcal{D}}^\Lambda$,

(13) $$\left\{ d\alpha_1 \oplus \cdots \oplus d\alpha_r \in \overline{\mathcal{D}} \mid (\operatorname{supp} d\alpha_1, \ldots, \operatorname{supp} d\alpha_r) = (S_1, \ldots, S_r) \right\},$$

which is isomorphic to an open disc of dimension

$$\sum_{n=1}^r |S_n| - 1 = d - r.$$

A member of (13) has the form

(14) $$d\alpha_1 \oplus \cdots \oplus d\alpha_r = \sum_{\lambda_n \in S_1} w_n \delta(\lambda - \lambda_n) \oplus \cdots \oplus \sum_{\lambda_n \in S_r} w_n \delta(\lambda - \lambda_n),$$

the weights $w_n$ of which serve as coordinates on (13). Indeed the notation is simplest if we fix an indexing of $\Lambda$,
$$\Lambda = \{\lambda_1, \ldots, \lambda_d\},$$
and then work with partitions of
$$[d] = \{1, \ldots, d\}$$
rather than of $\Lambda$. We set up notation for our coordinate space as follows.

Given a non-empty set $S \subseteq [d]$, we write $W^S$ for the projective space
$$W^S = \{w : S \to \mathbf{R} \mid \operatorname{supp} w \neq \emptyset\}/\sim,$$
where $w \sim w'$ if $w = rw'$ for some non-zero $r \in \mathbf{R}$. Thus $W^S$ consists of sequences indexed by $S$, distinguished up to proportionality, and is isomorphic to $\mathbf{RP}^{|S|-1}$. It is important to distinguish $W^S$ from $W^{S'}$ whenever $S \neq S'$. When it is desirable to emphasize the domain of $w : S \to \mathbf{R}$, we append it as a superscript, writing $w^S$. As usual, we write $w_n$, or $w_n^S$, to denote the value of $w$ at $n \in S$. We will always use homogeneous coordinates on $W^S$, treating $w \in W^S$ as a particular map $w : S \to \mathbf{R}$, and letting the equivalence relation $\sim$ be understood. We write $W_+^S$ for the open subset of $W^S$ defined by
$$W_+^S = \{w : S \to \mathbf{R} \mid w_m w_n > 0 \ \forall m, n \in S\}.$$

Referring to (13), let $(S_1', \ldots, S_r') \in \mathcal{P}_{[d]}$ be the partition corresponding to $(S_1, \ldots, S_r) \in \mathcal{P}_\Lambda$, so that for each $n$,
$$S_n' = \{i \in [d] \mid \lambda_i \in S_n\}.$$
Then the weights $w_n$ in (14) belong to the coordinate patch

(15) $$W_+^{(S_1', \ldots, S_r')} = W_+^{S_1'} \oplus \cdots \oplus W_+^{S_r'},$$

and we may identify (13) with $W_+^{(S_1', \ldots, S_r')}$. We denote an element
$$w^{S_1} \oplus \cdots \oplus w^{S_r} \in W_+^{(S_1, \ldots, S_r)},$$
using the abbreviated notation
$$w^{(S_1, \ldots, S_r)} = w^{S_1} \oplus \cdots \oplus w^{S_r}.$$

Now, $W_+^{(S_1', \ldots, S_r')}$ is a direct product of open subsets of finite-dimensional projective space, and carries its own topology. By the natural topology of (13), we mean the topology of this direct product. We identify $\overline{\mathcal{D}}^\Lambda$ with the union of coordinate patches (15), indexed by partitions of $[d]$,
$$\overline{\mathcal{D}}^\Lambda = \bigcup_{(S_1, \ldots, S_r) \in \mathcal{P}_{[d]}} W_+^{(S_1, \ldots, S_r)},$$

and use the obvious notation for the function $\overline{f} : \overline{\mathcal{D}}^\Lambda \to \overline{\mathcal{J}}^\Lambda$. Namely, for $w^S \in W_+^S$, $f(w^S)$ denotes

$$f\left(\sum_{n \in S} w_n^S \delta(\lambda - \lambda_n)\right),$$

and, given a point $w^{(S_1,\ldots,S_r)} \in W_+^{(S_1,\ldots,S_r)}$, we write

$$\overline{f}(w^{(S_1,\ldots,S_r)}) = f(w^{S_1}) \oplus \cdots \oplus f(w^{S_r}).$$

It is easy to see that the topology on $\overline{\mathcal{D}}^\Lambda$ induced by $\overline{f}$, restricted to a given piece $W_+^{(S_1,\ldots,S_r)}$, agrees with the natural topology. The restriction to $W_+^{(S_1,\ldots,S_r)}$ of each component $\overline{f}_n$ of $\overline{f}$ is either a rational function without singularities, or the square root of such a function, and the components $w_n$ of the normalized representor of the of the restriction of $\overline{f}^{-1}$ to $\overline{f}(W_+^{(S_1,\ldots,S_r)})$ are polynomials in the matrix entries $a_1, \ldots, a_{2d-1}$. Therefore the restriction to $W_+^{(S_1,\ldots,S_r)}$ of $\overline{f}$ is a homeomorphism. But there is no natural topology on $\overline{\mathcal{D}}^\Lambda$ as a whole; we need to analyze $\overline{f}$ directly to determine how its pieces fit together in the induced topology. Note that $\overline{\mathcal{D}}^\Lambda$ has a single piece, $W_+^{[d]}$, of maximum dimension $d-1$, corresponding to the trivial ordered partition $([d]) \in \mathcal{P}_{[d]}$. We refer to $W_+^{[d]}$ as the principal component of $\overline{\mathcal{D}}^\Lambda$. It turns out that the induced structure of $\overline{\mathcal{D}}^\Lambda$ is determined by sequences, or curves, in the principal component which tend to the boundary of $W_+^{[d]}$ in $W^{[d]}$. In fact it is enough just to consider moment curves in $W_+^{[d]}$ tending to the boundary. Note that $W_+^{[d]}$ may be thought of as the standard, open, $(d-1)$-dimensional simplex, as each point in $W_+^{[d]}$ has a unique representative $w$ satisfying $\sum_{n=1}^d w_n = 1$, with each $w_n > 0$.

2.3. **A blow-up of the standard simplex.** Given a $d$-element set

$$\Lambda = \{\lambda_1, \ldots, \lambda_d\} \subset \mathbf{R},$$

we construct a pair $(\mathcal{B}^\Lambda, \pi)$, consisting of a semi-algebraic set $\mathcal{B}^\Lambda$ and a bijective map

$$\pi : \mathcal{B}^\Lambda \to \overline{\mathcal{D}}^\Lambda.$$

The set $\mathcal{B}^\Lambda$ is defined in a relatively high-dimensional projective space $\mathcal{W}$ by quadratic equations and inequalities. The $\mathcal{W}$-structure of $\mathcal{B}^\Lambda$ encodes the topology on $\overline{\mathcal{D}}^\Lambda$ induced by $\overline{f}$ in that the map

$$\overline{f} \circ \pi : \mathcal{B}^\Lambda \to \overline{\mathcal{J}}^\Lambda$$

is a homeomorphism. (This latter fact will be proved in Section 2.7.) The pair $(\mathcal{B}^\Lambda, \pi)$ is a blow-up of $W_+^{(S_1,\ldots,S_r)} \subset \overline{\mathcal{D}}^\Lambda$ with centre

$$D^{(S_1,\ldots,S_r)} = \text{ the boundary of } W_+^{(S_1,\ldots,S_r)} \text{ in } W^{(S_1,\ldots,S_r)},$$

for every ordered partition $(S_1, \ldots, S_r) \in \mathcal{P}_{[d]}$, in the sense that there is a coordinate projection $\pi^{(S_1,\ldots,S_r)}$ on $\mathcal{W}$ such that:

- $\pi^{(S_1,\ldots,S_r)}(\mathcal{B}^\Lambda) = W_+^{(S_1,\ldots,S_r)} \cup D^{(S_1,\ldots,S_r)}$;

- the restriction of $\pi$ to $\pi^{-1}(W_+^{(S_1,\ldots,S_r)})$ coincides with $\pi^{(S_1,\ldots,S_r)}$.

In particular, $(\mathcal{B}^\Lambda, \pi)$ is a blow-up of the standard simplex $W_+^{[d]} \subset \overline{\mathcal{D}}^\Lambda$.

We now describe $(\mathcal{B}^\Lambda, \pi)$ in detail. $\mathcal{B}^\Lambda$ is contained in the ambient space

$$\mathcal{W} = \prod_{\emptyset \neq S \subseteq [d]} W^S.$$

The dimension of $\mathcal{W}$ is easily computed to be $(d-2)2^{d-1} + 1$. Given $w \in \mathcal{W}$, we write $w^S$ for its coordinate projection onto the factor $W^S$, and $w_n^S$ for the entry of $w^S$ corresponding to $n \in S$. Note that, as established earlier, these are homogeneous coordinates. The equations defining $\mathcal{B}^\Lambda$ involve the following scalars whose provenance we shall explain later, in Section 2.6. For each non-empty $S \subseteq [d]$ and each $i \in [d]$, set

$$C_i^S = \prod_{n \in S} (\lambda_n - \lambda_i)^2 \quad \text{and} \quad C_i^\emptyset = 1.$$

**Definition 2.** *Let $\mathcal{B}^\Lambda$ denote the semi-algebraic subset of $\mathcal{W}$ consisting of points $w$ which satisfy the equations and inequalities*

(16) $\qquad C_i^{S-R} w_j^R w_i^S = C_j^{S-R} w_i^R w_j^S \qquad \forall R \subseteq S \subseteq [d], \ i,j \in R$

(17) $\qquad w_i^S w_j^S \geq 0 \qquad\qquad\qquad\quad \forall S \subseteq [d], \ i,j \in S.$

We need some technical preparation in order to define the map $\pi : \mathcal{B}^\Lambda \to \overline{\mathcal{D}}^\Lambda$. Here it will be useful to consider the collection $\mathcal{L}_{[d]}$ of subsets of $[d]$, partially ordered by the inclusion relation $\subseteq$; we make use of a correspondence between totally ordered sets, or chains, in $\mathcal{L}_{[d]}$ having maximum element $[d]$, and ordered partitions of $[d]$. Let $w \in \mathcal{B}^\Lambda, S \subseteq [d], n \in S$, and write $P(S, n; w)$ for the statement

$$w_n^S \neq 0 \quad \text{and} \quad \forall T \subseteq [d], w_n^T \neq 0 \Rightarrow |T| \leq |S|.$$

**Lemma 3.** *If $P(R, m; w)$ and $P(S, n; w)$ then either $R \subseteq S$ or $S \subseteq R$.*

*Proof.* Set $T = R \cup S$ and choose $i \in T$ such that $w_i^T \neq 0$. (Such an $i$ exists because $w_n^T$ are homogeneous coordinates for the projective space $W^T$.) If $i \in R$, then $C_i^{T-R} w_m^R w_i^T \neq 0$, so by equation (16), $C_m^{T-R} w_i^R w_m^T \neq 0$, which implies $w_m^T \neq 0$. Therefore, since $P(R, m; w)$, $|T| \leq |R|$, forcing $T = R$, or equivalently, $S \subseteq R$. Similarly, if $i \in S$, then $P(S, n; w)$ and equation (16) force $R \subseteq S$. ∎

It follows from Lemma 3 that for each $n \in [d]$ there is a unique $S \subseteq [d]$ such that $P(S, n; w)$, so that we may define

$$\varphi_w : [d] \to \mathcal{L}_{[d]}$$

by setting

$$\varphi_w(n) = \text{ the unique } S \subseteq [d] \text{ such that } P(S, n; w).$$

Note that $[d]$ itself is in the range of $\varphi_w$, since necessarily $w_i^{[d]} \neq 0$, for some $i$, whence $P([d], i; w)$ automatically holds. Moreover, again by Lemma 3, the range $\{K_1, \ldots, K_r\} \subseteq \mathcal{L}_{[d]}$ of $\varphi_w$ is totally ordered; i.e., for an appropriate relabeling it is a chain of the form

$$[d] = K_1 \supset \cdots \supset K_r \neq \emptyset, \tag{18}$$

for some $r \geq 1$. We are now in a position to define the map $\pi$, but first we state another technical fact, for the sake of future reference.

**Lemma 4.** *Let $\{K_1, \ldots, K_r\} \subset \mathcal{L}_{[d]}$ be the range of $\varphi_w$, labeled so that (18) holds, and let $(S_1, \ldots, S_r) \in \mathcal{P}_{[d]}$ be the corresponding ordered partition, defined by*

$$(S_1, \ldots, S_r) = (K_1 - K_2, \ldots, K_{r-1} - K_r, K_r).$$

*Then for each $n$, $1 \leq n \leq r$,*

$$\varphi_w^{-1}(K_n) = S_n.$$

*Proof.* We claim that for each $i \in [d]$, $\varphi_w(i)$ is the smallest set $K_n$ containing $i$, a fact which implies the lemma. To prove the claim, suppose $i \in K_m \subseteq K_n$ and $\varphi_w(i) = K_n$, and choose $j \in [d]$ such that $\varphi_w(j) = K_m$; necessarily $j \in K_m$. Then

$$C_i^{K_n - K_m} w_j^{K_m} w_i^{K_n} \neq 0,$$

so by equation (16),

$$C_j^{K_n - K_m} w_i^{K_m} w_j^{K_n} \neq 0,$$

forcing $w_j^{K_n} \neq 0$. But this implies $|K_n| \leq |K_m|$, since by assumption $P(K_m, j; w)$ holds, forcing $K_n = K_m$. The claim follows. ∎

We define $\pi : \mathcal{B}^\Lambda \to \overline{\mathcal{D}}^\Lambda$ as follows.

**Definition 5.** *Let $w \in \mathcal{B}^\Lambda$, let $\{K_1, \ldots, K_r\}$ denote the range of $\varphi_w$, indexed such that (18) holds, and let*

$$(S_1, \ldots, S_r) = (K_1 - K_2, \ldots, K_{r-1} - K_r, K_r) \in \mathcal{P}_{[d]}$$

*be the corresponding ordered partition. Set*

$$\pi(w) = (w_n^{K_1})_{n \in S_1} \oplus \cdots \oplus (w_n^{K_r})_{n \in S_r} \in W_+^{(S_1, \ldots, S_r)}.$$

Note that the inequalities (17) in the definition of $\mathcal{B}^\Lambda$ guarantee that $\pi(w)$ is in $W_+^{(S_1, \ldots, S_r)}$ and not just in $W^{(S_1, \ldots, S_r)}$. Furthermore, because $(S_1, \ldots, S_r)$ determines the corresponding chain $K_1 \supset \cdots \supset K_r$ via

$$K_n = S_n \cup S_{n+1} \cup \ldots \cup S_r,$$

the restriction of $\pi$ to $\pi^{-1}(W_+^{(S_1, \ldots, S_r)})$ is a coordinate projection. Next we define a map $\rho : \overline{\mathcal{D}}^\Lambda \to \mathcal{B}^\Lambda$, which will turn out to be the inverse of $\pi$.

**Definition 6.** *Let $w^{(S_1, \ldots, S_r)} \in \overline{\mathcal{D}}^\Lambda$, and let $\{K_1, \ldots, K_r\}$ be the corresponding chain in $\mathcal{L}_{[d]}$, defined by*

$$K_n = S_n \cup S_{n+1} \cup \cdots \cup S_r$$

*for each $n$, $1 \leq n \leq r$. Set*

$$\rho(w^{(S_1, \ldots, S_r)}) = \widetilde{w},$$

where $\widetilde{w} \in \mathcal{W}$ is defined in terms of its components $\widetilde{w}^S$ as follows. Given $S \subseteq [d]$, choose $K_i \supseteq S$ to be the smallest member of the above chain that includes $S$, and set

$$\widetilde{w}_n^S = \begin{cases} C_n^{K_i-S} w_n^{S_i} & \text{if } n \in S_i \\ 0 & \text{if } n \notin S_i \end{cases}. \tag{19}$$

To verify that $\widetilde{w} = \rho(w^{(S_1,\ldots,S_r)}) \in \mathcal{B}^\Lambda$, we check that equations (16) and inequalities (17) of Definition 2 are satisfied. Indeed $\widetilde{w}$ satisfies the inequalities (17) by (19) and the definition of $W_+^{(S_1,\ldots,S_r)}$. To see that $\widetilde{w}$ satisfies equations (16), let $R \subseteq S \subseteq [d]$ with $i, j \in R$. Let $K_m \supseteq R$ be the smallest set in $\{K_1,\ldots,K_r\}$ that includes $R$, and let $K_n \supseteq S$ be the smallest that includes $S$; thus $K_m \subseteq K_n$. Applying Definition 6, the left-hand side of (16) works out to be

$$C_i^{S-R} \widetilde{w}_j^R \widetilde{w}_i^S = \begin{cases} C_i^{S-R}(C_j^{K_m-R} w_j^{S_m})(C_i^{K_n-S} w_i^{S_n}) & \text{if } j \in S_m \text{ and } i \in S_n \\ 0 & \text{otherwise} \end{cases}$$

$$= \begin{cases} C_j^{K_m-R} C_i^{K_n-R} w_j^{S_m} w_i^{S_n} & \text{if } j \in S_m \text{ and } i \in S_n \\ 0 & \text{otherwise} \end{cases}, \tag{20}$$

the latter equality since $C_i^{K_n-S} C_i^{S-R} = C_i^{K_n-R}$. The right-hand side of (16) is

$$C_j^{S-R} \widetilde{w}_i^R \widetilde{w}_j^S = \begin{cases} C_j^{S-R}(C_i^{K_m-R} w_i^{S_m})(C_j^{K_n-S} w_j^{S_n}) & \text{if } i \in S_m \text{ and } j \in S_n \\ 0 & \text{otherwise} \end{cases}$$

$$= \begin{cases} C_j^{K_n-R} C_i^{K_m-R} w_i^{S_m} w_j^{S_n} & \text{if } i \in S_m \text{ and } j \in S_n \\ 0 & \text{otherwise} \end{cases}, \tag{21}$$

since $C_j^{K_n-S} C_j^{S-R} = C_j^{K_n-R}$. Now, if $m \neq n$, then $S_n \subseteq K_n - K_m$, and therefore, since $i, j \in R \subseteq K_m$, neither $i$ nor $j$ belongs to $S_n$ and consequently both (20) and (21) are 0. On the other hand, if $m = n$, then (20) and (21) are precisely the same. So equations (16) are satisfied in any case, confirming that $\rho(w^{(S_1,\ldots,S_r)}) \in \mathcal{B}^\Lambda$.

**Proposition 7.** *The maps $\rho$ and $\pi$ are inverses, that is, $\rho \circ \pi = 1_{\mathcal{B}^\Lambda}$ and $\pi \circ \rho = 1_{\overline{\mathcal{D}}^\Lambda}$.*

*Proof.* ($\pi \circ \rho = 1_{\overline{\mathcal{D}}^\Lambda}$.) Let $w^{(S_1,\ldots,S_r)} \in \overline{\mathcal{D}}^\Lambda$ and set $\widetilde{w} = \rho(w^{(S_1,\ldots,S_r)})$; we need to check that $\pi(\widetilde{w}) = w^{(S_1,\ldots,S_r)}$. It follows directly from (19) that for $n \in S_i$,

$$\varphi_{\widetilde{w}}(n) = K_i = S_i \cup S_{i+1} \cup \ldots \cup S_r.$$

Then, by Definition 5,

$$\pi(\widetilde{w}) = (\widetilde{w}_n^{K_1})_{n \in S_1} \oplus \cdots \oplus (\widetilde{w}_n^{K_r})_{n \in S_r},$$

where, again by (19), $\widetilde{w}_n^{K_i} = C_n^{K_i-K_i} w_n^{S_i} = w_n^{S_i}$, so that $\pi(\widetilde{w}) = w^{(S_1,\ldots,S_r)}$.

($\rho \circ \pi = 1_{\mathcal{B}^\Lambda}$.) Let $\widetilde{w} \in \mathcal{B}^\Lambda$ and let $w^{(S_1,\ldots,S_r)} = \pi(\widetilde{w})$ be its image by $\pi$, so that $w_n^{S_j} = \widetilde{w}_n^{K_j}$, where $K_j = S_j \cup S_{j+1} \cup \ldots \cup S_r$. In order to verify that $\rho(w^{(S_1,\ldots,S_r)}) = \widetilde{w}$, fix $S \subseteq [d]$, and, following Definition 6, choose

$$K_i = S_i \cup S_{i+1} \cup \ldots \cup S_r \supseteq S$$

to be the smallest such set that includes $S$, so that $S \cap S_i \neq \emptyset$. Then

$$\rho(w^{(S_1,\ldots,S_r)})_n^S = \begin{cases} C_n^{K_i-S} w_n^{S_i} & \text{if } n \in S_i \\ 0 & \text{if } n \notin S_i \end{cases}$$

(22)
$$= \begin{cases} C_n^{K_i-S} \widetilde{w}_n^{K_i} & \text{if } n \in S_i \\ 0 & \text{if } n \notin S_i \end{cases}.$$

For comparison, fix $m \in S \cap S_i$, and compute $\widetilde{w}_n^S$ using equation (16):

$$C_m^{K_i-S} \widetilde{w}_n^S \widetilde{w}_m^{K_i} = C_n^{K_i-S} \widetilde{w}_m^S \widetilde{w}_n^{K_i}$$

(23)
$$\Rightarrow \quad \widetilde{w}_n^S = \left( \frac{\widetilde{w}_m^S}{C_m^{K_i-S} \widetilde{w}_m^{K_i}} \right) C_n^{K_i-S} \widetilde{w}_n^{K_i}.$$

We will verify that the representations (22) and (23) are proportional, which implies that $\rho(w^{(S_1,\ldots,S_r)})^S = \widetilde{w}^S$, as desired. By Lemma 4, if $n \notin S_i$, that is, if $n \in S_{i+1} \cup \ldots \cup S_r$, then

$$\varphi_{\widetilde{w}}(n) \in \{K_{i+1}, \ldots, K_r\}$$

and consequently $\widetilde{w}_n^{K_i} = 0$. On the other hand, if $n \in S_i$, then $\varphi_{\widetilde{w}}(n) = K_i$ and $\widetilde{w}_n^{K_i} \neq 0$; in particular, $\widetilde{w}_m^{K_i} \neq 0$. The coefficient $\widetilde{w}_m^S/(C_m^{K_i-S} \widetilde{w}_m^{K_i})$ in (23) is independent of $n \in S$, and non-zero since by definition $\widetilde{w}^S$ has at least one non-zero entry, therefore (22) and (23) are proportional. ∎

2.4. **Stable sequences.** We focus our analysis of $\overline{f} : \overline{\mathcal{D}}^\Lambda \to \overline{\mathcal{J}}^\Lambda$ mainly on its values at sequences in $W_+^{[d]}$, the principal component of $\overline{\mathcal{D}}^\Lambda$, that tend to the boundary of $W_+^{[d]}$ in $W^{[d]}$. The blow-up $(\mathcal{B}^\Lambda, \pi)$ of $W_+^{[d]}$ has a central place in this analysis; indeed, as mentioned earlier, we are working toward a proof that

$$\overline{f} \circ \pi : \mathcal{B}^\Lambda \to \overline{\mathcal{J}}^\Lambda$$

is a homeomorphism. We keep the notation $\rho = \pi^{-1}$, consistent with Definition 6 and Proposition 7. Let $\mathcal{B}_\circ^\Lambda$ denote the piece of $\mathcal{B}^\Lambda$ corresponding to $W_+^{[d]}$, that is,

$$\mathcal{B}_\circ^\Lambda = \rho(W_+^{[d]}).$$

Note that $\mathcal{B}_\circ^\Lambda$ is characterized as the collection of points $w \in \mathcal{B}^\Lambda$ for which each entry $w_n^S$, $n \in S \subseteq [d]$, is non-zero. For $w \in \mathcal{B}_\circ^\Lambda$ the range of $\varphi_w$ is the singleton set $\{[d]\}$. The restriction of $\pi$ to $\mathcal{B}_\circ^\Lambda$ coincides with projection of $w$ onto the coordinates $w^{[d]}$, and is a homeomorphism with respect to the natural topology of $W_+^{[d]}$.

**Definition 8.** *Call a sequence $w(n) \in W_+^{[d]}$, $n \in \mathbf{Z}^+$, stable if the corresponding sequence $\rho(w(n)) \in \mathcal{B}_\circ^\Lambda$ is convergent in $\mathcal{W}$ (equivalently, in $\mathcal{B}^\Lambda$).*

By Definition 6, convergence of $\rho(w(n))$ in $\mathcal{W}$ amounts to the simultaneous convergence of the sequences

$$\left( C_i^{[d]-S} w(n)_i \right)_{i \in S} \quad \text{in } W^S,$$

for each non-empty $S \subseteq [d]$. Since the constants $C_i^{[d]-S}$ ($i \in S$) are strictly positive, this is equivalent to the simultaneous convergence of
$$(w(n)_i)_{i \in S} \text{ in } W^S.$$
Thus we have the following characterization. A sequence $w(n) \in W_+^{[d]}$ is stable if and only if, for every $i, j \in [d]$, there exists $L_{ij} \in \mathbf{R}^+ \cup \{0, \infty\}$ such that
$$\frac{w(n)_i}{w(n)_j} \to L_{ij} \text{ as } n \to \infty.$$
The $L_{ij}$ determine an ordered partition $(S_1, \ldots, S_r) \in \mathcal{P}_{[d]}$, where cells of the partition are defined by the equivalence relation
$$i \sim j \Leftrightarrow 0 < L_{ij} < \infty,$$
and the ordering of the cells $S_m < S_n$ is determined by
$$i \in S_m, j \in S_n \Rightarrow L_{ij} = \infty \text{ (equivalently } L_{ji} = 0).$$
We call $(S_1, \ldots, S_r)$ the limiting partition of the stable sequence $w(n) \in W_+^{[d]}$; it is also characterized by the fact that
$$\pi\left(\lim_{n \to \infty} \rho(w(n))\right) \in W_+^{(S_1, \ldots, S_r)}.$$

Let us use the notation
$$\bigoplus_{(S_1, \ldots, S_r)} w^{S_i} = w^{S_1} \oplus \cdots \oplus w^{S_r}$$
for $w^{S_1} \oplus \cdots \oplus w^{S_r} \in W^{(S_1, \ldots, S_r)}$. Given a stable sequence $w(n)$ in $W_+^{[d]}$ with limiting distribution $(S_1, \ldots, S_r)$, we calculate
$$\pi\left(\lim_{n \to \infty} \rho(w(n))\right)$$
as follows. As usual, for each $i$, $1 \leq i \leq r$, set
$$K_i = S_i \cup S_{i+1} \cup \ldots S_r.$$
By Definition 5, the restriction of $\pi$ to $\rho\left(W_+^{(S_1, \ldots, S_r)}\right)$ coincides with a coordinate projection on $\mathcal{W}$, call it $\pi^{(S_1, \ldots, S_r)}$. Then

$$
\begin{aligned}
\pi\left(\lim_{n \to \infty} \rho(w(n))\right) &= \pi^{(S_1, \ldots, S_r)}\left(\lim_{n \to \infty} \rho(w(n))\right) \\
&= \lim_{n \to \infty} \pi^{(S_1, \ldots, S_r)}(\rho(w(n))) \\
&= \lim_{n \to \infty} \bigoplus_{(S_1, \ldots, S_r)} \left(\rho(w(n))_j^{K_i}\right)_{j \in S_i} \\
&= \lim_{n \to \infty} \bigoplus_{(S_1, \ldots, S_r)} \left(C_j^{[d]-K_i} w(n)_j\right)_{j \in S_i} \\
(24) \qquad &= \bigoplus_{(S_1, \ldots, S_r)} \lim_{n \to \infty} \left(C_j^{[d]-K_i} w(n)_j\right)_{j \in S_i}.
\end{aligned}
$$

Returning to distributional notation, write

$$d\alpha_i(n) = \sum_{j \in S_i} w(n)_j \, \delta(\lambda - \lambda_j),$$

so that $\operatorname{supp} d\alpha_i(n)$ is indexed by $S_i$, and let

$$p_{d\alpha_i}(\lambda) = \prod_{j \in S_i}(\lambda - \lambda_j)$$

denote the characteristic polynomial of $f(d\alpha_i(n))$, which does not depend on $n$. The assumption that $w(n)$ is stable guarantees that each

$$d\overline{\alpha}_i = \lim_{n \to \infty} d\alpha_i(n)$$

exists and has support indexed by $S_i$. Rewriting (24), we have

$$(25) \quad \pi\left(\lim_{n \to \infty} \rho(d\alpha_1(n) + \ldots + d\alpha_r(n))\right) = d\overline{\alpha}_1 \oplus p_{d\alpha_1}^2 \, d\overline{\alpha}_2 \oplus \cdots \oplus p_{d\alpha_1 + \ldots + d\alpha_{r-1}}^2 \, d\overline{\alpha}_r.$$

Now, the sequence

$$(26) \qquad d\alpha(n) = d\overline{\alpha}_1 + n^{-1} \, d\overline{\alpha}_2 + \ldots + n^{-r+1} \, d\overline{\alpha}_r$$

is stable, has limiting distribution $(S_1, \ldots, S_r)$, and $\rho(d\alpha(n))$ converges to the image of (25) by $\rho$. Replacing $n^{-1}$ in (26) with the continuous variable $t$, $0 < t < 1$, yields the following.

**Proposition 9.** *For every stable sequence $w(n)$ in $W_+^{[d]}$ there exists a corresponding moment curve in $W_+^{[d]}$ of the form*

$$d\alpha(t) = d\alpha_1 + t \, d\alpha_2 + \ldots + t^{r-1} \, d\alpha_r, \qquad 0 < t < 1,$$

*where $(\operatorname{supp} d\alpha_1, \ldots, \operatorname{supp} d\alpha_r) \in \mathcal{P}_\Lambda$, such that*

$$\begin{aligned} \pi\left(\lim_{n \to \infty} \rho(w(n))\right) &= \pi\left(\lim_{t \to 0^+} \rho(d\alpha(t))\right) \\ &= d\alpha_1 \oplus p_{d\alpha_1}^2 \, d\alpha_2 \oplus \cdots \oplus p_{d\alpha_1 + \ldots + d\alpha_{r-1}}^2 \, d\alpha_r. \end{aligned}$$

Thus in order to study the way the pieces of $\overline{\mathcal{D}}^\Lambda$ are embedded in $\mathcal{B}^\Lambda$ by $\rho$, it is enough to study the behaviour of moment curves. In particular,

**Proposition 10.** $\mathcal{B}^\Lambda = \operatorname{cl} \mathcal{B}_\circ^\Lambda$.

*Proof.* $\mathcal{B}^\Lambda$ is closed, so $\operatorname{cl} \mathcal{B}_\circ^\Lambda \subseteq \mathcal{B}^\Lambda$. To see that the reverse inclusion holds, let $w \in \mathcal{B}^\Lambda$ and write

$$d\alpha_1 \oplus \cdots \oplus d\alpha_r = \pi(w).$$

Then $w$ is a limit point of $\mathcal{B}_\circ^\Lambda$, as we can see by choosing the appropriate moment curve:

$$w = \lim_{t \to 0^+} \rho\left(d\alpha_1 + t p_{d\alpha_1}^{-2} \, d\alpha_2 + \ldots + t^{r-1} p_{d\alpha_1 + \ldots + d\alpha_{r-1}}^{-2} \, d\alpha_r\right).$$

∎

Essentially the same argument as in the above proof, carried out within the lower-dimensional pieces $W_+^{(S_1,\ldots,S_r)}$ of $\overline{\mathcal{D}}^\Lambda$, reveals that $\mathcal{B}^\Lambda$ has a natural "facial structure" identical with that of a special convex polyhedron known as a permutahedron. The details of this combinatorial analysis appear in Section 3.

2.5. **Two technical lemmas.** In this section we derive two technical lemmas concerning stable sequences in $W_+^{[d]}$ by directly analyzing the maps

$$f : \mathcal{D} \to \mathcal{J}, \quad \Psi : \mathcal{D} \to MOP$$

of Definition 1.

To begin, we recall some necessary facts and set up notation. Given a matrix $J$, let $\chi(J)$ denote its characteristic polynomial, and let $J_i$ denote the leading $i \times i$ submatrix. Note that for $J = f(w)$, $w \in W_+^{[d]}$, the sequence $(p_0, \ldots, p_d) = \Psi(w)$ satisfies $p_i = \chi(J_i)$, $1 \leq i \leq d$. Given $w \in W_+^{[d]}$ and a non-empty set $Q \subset [d]$, let $w|_Q$ denote the restriction of $w : [d] \to \mathbf{R}$ to $Q$. Recall that, with the notation we have established earlier,

$$f(w|_Q) = f\left(\sum_{n \in Q} w_n \delta(\lambda - \lambda_n)\right).$$

The following notion will be used only fleetingly in the present section; in Section 3, however, it will be used extensively. A refinement of an ordered partition $(S_1, \ldots, S_r) \in \mathcal{P}_{[d]}$ is an ordered partition of the form

$$(R_1^1, \ldots, R_{i_1}^1, R_1^2, \ldots, R_{i_2}^2, \ldots, R_1^r, \ldots, R_{i_r}^r),$$

where, for each $j$, $1 \leq j \leq r$, $(R_1^j, \ldots, R_{i_j}^j)$ is an ordered partition of $S_j$. We write

$$(R_1, \ldots, R_q) \preceq (S_1, \ldots, S_r)$$

to indicate that $(R_1, \ldots, R_q)$ is a refinement of $(S_1, \ldots, S_r)$; conversely, we say that $(S_1, \ldots, S_r)$ is a coarsening of $(R_1, \ldots, R_q)$.

We note a pair of bounds on the values of $\Psi$ on $W_+^{[d]}$, in preparation for the coming lemmas. Let $P, S \subset \mathbf{R}$ denote the sets

$$P = \left\{ p_i^2(\lambda_k) \mid p_i = \Psi_i(w), \ w \in W_+^{[d]}, \ \lambda_k \in \Lambda, \ 0 \leq i \leq d \right\},$$

$$S = \left\{ \sum_{k \in Q} p_i^2(\lambda_k) \mid p_i = \Psi_i(w), \ w \in W_+^{[d]}, \ Q \subseteq [d], \ i < |Q| \right\}.$$

**Proposition 11.** *Let $U = \sup P$ and $L = \inf S$. Then $U < \infty$ and $L > 0$.*

*Proof.* Let $S' \subset \mathbf{R}$ denote the set

$$S' = \left\{ \sum_{k \in Q} p_i^2(\lambda_k) \mid p_i = \chi(J_i), \ J \in \overline{\mathcal{J}}^\Lambda, \ Q \subseteq [d], \ i < |Q| \right\},$$

and observe that $S' \supseteq S$. Now, in $S'$, each $\sum_{k \in Q} p_i^2(\lambda_k) > 0$, since $p_i$ has degree $i$ and $|Q| > i$. Moreover, the map $\chi(J)$ is continuous in the entries of $J$, and the domain $\overline{\mathcal{J}}^\Lambda$ is compact, unlike $W_+^{[d]}$. Therefore

$$L' = \inf S' > 0,$$

which implies positivity of $L \geq L'$, since $S \subseteq S'$. That $U < \infty$ also follows from compactness of $\overline{\mathcal{J}}^\Lambda$, by a similar argument. ∎

Note that if $(S_1, \ldots, S_r) \in \mathcal{P}_{[d]}$ is the limiting partition of a stable sequence $w(n) \in W_+^{[d]}$, then for any coarsening $(Q, R) \in \mathcal{P}_{[d]}$ of $(S_1, \ldots, S_r)$,

$$\lim_{n \to \infty} \frac{w(n)_j}{w(n)_i} = 0 \quad \forall \, i \in Q, \, j \in R.$$

Of course, $r \geq 2$ is necessary for such a coarsening to exist.

**Lemma 12.** *Let $w(n)$ be a stable sequence in $W_+^{[d]}$, with limiting partition $(S_1, \ldots, S_r)$, where $r \geq 2$, and let $(Q, R) \in \mathcal{P}_{[d]}$ be a coarsening of $(S_1, \ldots, S_r)$. If $f(w(n)|_Q)$ converges, then, for each $N$, $1 \leq N \leq 2|Q| - 1$, the sequence $f_N(w(n))$ converges, and*

$$\tag{27} \lim_{n \to \infty} f_N(w(n)) = \lim_{n \to \infty} f_N(w(n)|_Q).$$

*Proof.* We proceed by induction on $N$, in the range $1 \leq N \leq 2|Q| - 1$. We first consider the case $N = 1$. By Definition 1,

$$f_1(w(n)) = \frac{A + B}{C + D}, \quad f_1(w(n)|_Q) = \frac{A}{C},$$

where

$$A = \sum_{i \in Q} \lambda_i w(n)_i, \quad B = \sum_{i \in R} \lambda_i w(n)_i,$$

$$C = \sum_{i \in Q} w(n)_i, \quad D = \sum_{i \in R} w(n)_i.$$

In order to show that

$$\lim_{n \to \infty} f_1(w(n)) = \lim_{n \to \infty} f_1(w(n)|_Q),$$

it suffices to show that

$$\tag{28} \left| \frac{A+B}{C+D} - \frac{A}{C} \right| = \left| \frac{B}{C+D} - \frac{A}{C}\left(\frac{D}{C+D}\right) \right| \to 0,$$

as $n \to \infty$. Now, $A/C = f_1(w(n)|_Q)$ converges in $\mathbf{R}$ by assumption; and we know that for every $i \in Q$, $j \in R$, $\lim_{n \to \infty} w(n)_j/w(n)_i = 0$, which implies

$$\lim_{n \to \infty} \frac{B}{C+D} = \lim_{n \to \infty} \frac{D}{C+D} = 0,$$

and hence (28).

Next, suppose that (27) holds for every $N$ in the range $1 \leq N \leq 2i - 3$, where $2 \leq i \leq |Q|$; we will prove that (27) then holds for $N = 2i - 2$ and $N = 2i - 1$ also. In the case $N = 2i - 2$, it is simpler to prove the equivalent statement to (27), that

$$\text{(29)} \qquad \lim_{n \to \infty} f_{2i-2}^2(w(n)) = \lim_{n \to \infty} f_{2i-2}^2(w(n)|_Q).$$

Write

$$(p_0^n, \ldots, p_d^n) = \Psi(w(n)), \quad (q_0^n, \ldots, q_{|Q|}^n) = \Psi(w(n)|_Q),$$

where here the superscript $n$ in $p_i^n$ and $q_i^n$ is not an exponent, but just an index. By Definition 1,

$$f_{2i-2}^2(w(n)) = \frac{A + B}{C + D}, \quad f_{2i-2}^2(w(n)|_Q) = \frac{A'}{C'},$$

where

$$A = \sum_{k \in Q} (p_{i-1}^n(\lambda_k))^2 w(n)_k, \quad B = \sum_{k \in R} (p_{i-1}^n(\lambda_k))^2 w(n)_k,$$

$$C = \sum_{k \in Q} (p_{i-2}^n(\lambda_k))^2 w(n)_k, \quad D = \sum_{k \in R} (p_{i-2}^n(\lambda_k))^2 w(n)_k,$$

$$A' = \sum_{k \in Q} (q_{i-1}^n(\lambda_k))^2 w(n)_k, \quad C' = \sum_{k \in Q} (q_{i-2}^n(\lambda_k))^2 w(n)_k.$$

In order to prove (29), it suffices to show that the right-hand side of the inequality

$$\text{(30)} \qquad \left| \frac{A + B}{C + D} - \frac{A'}{C'} \right| \leq \left| \frac{A + B}{C + D} - \frac{A}{C} \right| + \left| \frac{A}{C} - \frac{A'}{C'} \right|$$

tends to 0 as $n \to \infty$.

Since $p_{i-1}^n = \chi(J_{i-1})$, where $J = f(w(n))$, and $q_{i-1}^n = \chi(J'_{i-1})$, where $J' = f(w(n)|_Q)$, the induction hypothesis and continuity of $\chi$ imply that $|A/A'| \to 1$ as $n \to \infty$. If $i = 2$, then trivially $p_{i-2}^n = q_{i-2}^n$; if $i > 2$, the same reasoning as above shows that $|C/C'| \to 1$ as $n \to \infty$. Since $A'/C'$ converges by assumption, it follows that

$$\lim_{n \to \infty} \left| \frac{A}{C} - \frac{A'}{C'} \right| = 0,$$

and hence that $A/C$ converges. We now argue that the remaining term on the right-hand side of (30),

$$\text{(31)} \qquad \left| \frac{A + B}{C + D} - \frac{A}{C} \right| = \left| \frac{B}{C + D} - \frac{A}{C} \left( \frac{D}{C + D} \right) \right|,$$

tends to 0 as $n \to \infty$. For each $n$, let $w(n)_{\min} = w(n)_{k(n)}$, where the index $k(n) \in Q$ renders $|w(n)_{k(n)}|$ minimal over $Q$ for the given $n$. Then, by Proposition 11,

$$B/w(n)_{\min} \leq U \sum_{k \in R} w(n)_k / w(n)_{\min} \to 0 \quad \text{as } n \to \infty,$$

$$D/w(n)_{\min} \leq U \sum_{k \in R} w(n)_k / w(n)_{\min} \to 0 \quad \text{as } n \to \infty,$$

$$C/w(n)_{\min} \geq \quad L \quad > 0 \quad \text{independently of } n.$$

It follows that (31) tends to 0 as $n \to \infty$, proving (29). We omit the proof that

$$\lim_{n \to \infty} f_{2i-1}(w(n)) = \lim_{n \to \infty} f_{2i-1}(w(n)|_Q),$$

which goes along exactly the same lines as the proof of (29), with obvious minor modifications. ∎

**Lemma 13.** *Let $w(n)$ be a stable sequence in $W_+^{[d]}$, with limiting partition $(S_1, \ldots, S_r)$, where $r \geq 2$. Then*

$$\lim_{n \to \infty} f_{2|S_1|}(w(n)) = 0.$$

*Proof.* To simplify the notation, write $N = |S_1|$, and set $Q = S_1$ and $R = S_2 \cup \ldots \cup S_r$. As before, let $w(n)_{\min} = w(n)_{k(n)}$, where $k(n) \in Q$ renders $|w(n)_{\min}|$ minimal over $Q$ for the given $n$. The fact that $S_1$ is a cell in the limiting partition of $w(n)$ implies that

(32) $\qquad \forall k \in Q$, the sequence $w(n)_k / w(n)_{\min}$ is bounded

(indeed it is convergent). Set $p_N^n = \Psi_N(w(n))$ and $q_N^n = \Psi_N(w(n)|_Q)$, where the superscript $n$ in $p_N^n$ and $q_N^n$ is an index and not an exponent, so that $p_N^n = \chi(J_N)$, where $J = f(w(n))$, and $q_N^n = \chi(J')$, where $J' = f(w(n)|_Q)$. In fact $q_N^n$ has constant value

$$q_N^n(\lambda) = \prod_{k \in Q} (\lambda - \lambda_k),$$

independent of $n$. The entries of $J_N$, $J'$ are, respectively,

$$(f_1(w(n)), \ldots, f_{2N-1}(w(n))), \quad (f_1(w(n)|_Q), \ldots, f_{2N-1}(w(n)|_Q)).$$

Therefore, by Lemma 12, $\lim_{n \to \infty} J_N = \lim_{n \to \infty} J'$, and by continuity of $\chi$,

$$\lim_{n \to \infty} p_N^n(\lambda) = \lim_{n \to \infty} q_N^n(\lambda) = \prod_{k \in Q} (\lambda - \lambda_k).$$

Using (32), this implies

(33) $$\lim_{n \to \infty} \sum_{k \in Q} (p_N^n(\lambda_k))^2 w(n)_k / w(n)_{\min} = 0.$$

For every $k \in R$, $w(n)_k/w(n)_{\min} \to 0$, so, by Proposition 11,

$$\sum_{k \in R}(p_N^n(\lambda_k))^2 w(n)_k/w(n)_{\min} \leq U \sum_{k \in R} w(n)_k/w_{\min} \to 0, \tag{34}$$

as $n \to \infty$. On the other hand, since $N - 1 < |Q|$, we have, by Proposition 11,

$$\sum_{k \in Q}(p_{N-1}^n(\lambda_k))^2 w(n)_k/w(n)_{\min} \geq L > 0, \tag{35}$$

independently of $n$, where $p_{N-1}^n = \Psi_{N-1}(w(n))$. By Definition 1,

$$f_{2N}^2(w(n)) = \frac{\sum_{k \in Q}(p_N^n(\lambda_k))^2 w(n)_k + \sum_{k \in R}(p_N^n(\lambda_k))^2 w(n)_k}{\sum_{k \in Q}(p_{N-1}^n(\lambda_k))^2 w(n)_k + \sum_{k \in R}(p_{N-1}^n(\lambda_k))^2 w(n)_k}.$$

It follows from (33), (34) and (35) that $\lim_{n \to \infty} f_{2N}^2(w(n)) = 0$. ∎

2.6. **The flip transpose.** In this section we explain the origin of the coefficients $C_i^{S-R}$ appearing in the equations for $\mathcal{B}^\Lambda$ in Definition 2. The underlying technical fact is a formula involving the flip transpose.

The flip transpose of a $d \times d$ matrix $M$ is defined to be the $d \times d$ matrix $M^F$ having entries

$$M^F(i,j) = M(d+1-j, d+1-i).$$

In terms of our convention (1) for labeling the entries $(a_1, \ldots, a_{2d-1})$ of a Jacobi matrix $J$, the entries of $J^F$ are simply the reverse sequence $(a_{2d-1}, \ldots, a_1)$. The flip transpose of a matrix induces, in the obvious way, a corresponding operation on probability distributions in $\mathcal{D}$, namely, define the flip transpose of $d\alpha \in \mathcal{D}$, written $d\alpha^F$, to be the spectral distribution of

$$(f(d\alpha))^F.$$

Note that flip transposition leaves the spectrum of a matrix unchanged, so that $\operatorname{supp} d\alpha^F = \operatorname{supp} d\alpha$. As an operation on distributions, therefore, flip transposition carries $W_+^{[d]}$ into itself; we write $w^F$ for the flip transpose of $w \in W_+^{[d]}$. It follows immediately from the definition that for each $i$, $1 \leq i \leq 2d-1$,

$$f_i(w) = f_{2d-i}(w^F). \tag{36}$$

Define constants $\overline{C}_i^S$ in terms of the $d$-element set $\Lambda = \{\lambda_1, \ldots, \lambda_d\}$ as follows. Set

$$p(\lambda) = \prod_{n=1}^{d}(\lambda - \lambda_n), \quad p'(\lambda) = \frac{d}{d\lambda}p(\lambda),$$

and for each non-empty $S \subseteq [d]$, $i \in [d]$, set

$$\overline{C}_i^S = \prod_{n \in S}(p'(\lambda_i - \lambda_n))^2.$$

Then, for $i \in R \subseteq S \subseteq [d]$, we have the relation $\overline{C}_i^S/\overline{C}_i^R = C_i^{S-R}$.

The following result is a translation of [9, Proposition 3.3] into our current notation.

**Proposition 14.** *For $w \in W_+^{[d]}$, $w^F$ is given by the formula*

$$w_n^F = \frac{1}{\overline{C}_n^{[d]} w_n}, \quad \forall n \in [d].$$

This has an immediate implication for stable sequences.

**Corollary 15.** *Let $w(n)$ be a stable sequence in $W_+^{[d]}$ with limiting partition $(S_1, \ldots, S_r)$. Then the sequence $w(n)^F \in W_+^{[d]}$ is stable and has the reverse-ordered limiting partition $(S_r, \ldots, S_1)$.*

Let $Q \cup R$ be a partition of $[d]$, and let $w \in W_+^{[d]}$. Repeated application of Proposition 14 yields

$$\begin{aligned}
(w^F|_R)^F &= \left(\left(\frac{1}{\overline{C}_i^{[d]} w_i}\right)_{i \in R}\right)^F \\
&= \left(\frac{\overline{C}_i^{[d]} w_i}{\overline{C}_i^R}\right)_{i \in R} \\
&= (C_i^{[d]-R} w_i)_{i \in R} \\
&= (C_i^Q)_{i \in R}.
\end{aligned}$$
(37)

The following technical result accounts for the coefficients that appear in the equations defining $\mathcal{B}^\Lambda$.

**Lemma 16.** *Let $w(n)$ be a stable sequence in $W_+^{[d]}$ with limiting partition $(S_1, \ldots, S_r)$, where $r \geq 2$, and let $(Q, R) \in \mathcal{P}_{[d]}$ be a coarsening of $(S_1, \ldots, S_r)$. If $f(w(n)^F|_R)$ converges, then for each $i$, $2|Q| + 1 \leq i \leq 2d - 1$, $f_i(w(n))$ converges and*

$$\begin{aligned}
\lim_{n \to \infty} f_i(w(n)) &= \lim_{n \to \infty} f_{i-2|Q|}\left((w(n)^F|_R)^F\right) \\
&= \lim_{n \to \infty} f_{i-2|Q|}\left((C_j^Q w(n)_j)_{j \in R}\right).
\end{aligned}$$

*Proof.* By Corollary 15 $w(n)^F$ is stable with limiting partition $(S_r, \ldots, S_1)$, of which $(R, Q)$ is a coarsening. If $f(w(n)^F|_R)$ converges then, by Lemma 12, for each $i$, $1 \leq i \leq 2|R| - 1$, $f_i(w(n)^F)$ converges, and

$$\lim_{n \to \infty} f_i(w(n)^F) = \lim_{n \to \infty} f_i(w(n)^F|_R).$$

Applying (36) to the left- and right-hand sides yields

$$\lim_{n \to \infty} f_j(w(n)) = \lim_{n \to \infty} f_{j-2|Q|}((w(n)^F|_R)^F),$$

where $j = 2d - i$ and $2|Q| + 1 \leq j \leq 2d - 1$. The calculation (37) then gives the last expression in the conclusion of the lemma. ∎

**Lemma 17.** *Let $w(n)$ be a stable sequence in $W_+^{[d]}$ with limiting partition $(S_1, \ldots, S_r)$, where $r \geq 2$. Then*

$$\lim_{n \to \infty} f_{2d-2|S_r|}(w(n)) = 0.$$

*Proof.* Applying Lemma 13 to the sequence $w(n)^F$ yields

$$\lim_{n \to \infty} f_{2|S_r|}(w(n)^F) = 0.$$

By (36) this is equivalent to the conclusion of the present lemma. ∎

2.7. **The main result.** In this section we prove that

$$\overline{f} \circ \pi : \mathcal{B}^\Lambda \to \overline{\mathcal{J}}^\Lambda$$

is a homeomorphism.

The proof involves stable sequences in $W_+^Q$, where $Q$ ranges over non-empty subsets of $[d]$. We use the characterization of stable sequences in terms of limits, to avoid dealing with semi-algebraic sets of the form $\mathcal{B}^\Gamma$ where $\Gamma \subset \Lambda$. Thus a sequence $w(n) \in W_+^Q$ is stable if, for every $i, j \in Q$, there exists $L_{i,j} \in \mathbf{R}^+ \cup \{0, \infty\}$ such that

$$\frac{w(n)_i}{w(n)_j} \to L_{ij} \quad \text{as} \quad n \to \infty.$$

Note that if $w(n)$ is a stable sequence in $W_+^{[d]}$ then for $Q \subseteq [d]$, the restriction $w(n)|_Q$ is automatically a stable sequence in $W_+^Q$.

**Proposition 18.** *Let $w(n)$ be a stable sequence in $W_+^{[d]}$ with limiting partition $(S_1, \ldots, S_r)$. Then $f(w(n))$ converges, and*

$$(38) \qquad \lim_{n \to \infty} f(w(n)) = \overline{f}\left( \bigoplus_{(S_1, \ldots, S_r)} \lim_{n \to \infty} \left( C_i^{S_1 \cup \ldots \cup S_{j-1}} w(n)_i \right)_{i \in S_j} \right).$$

*Proof.* We argue by induction on $N = |S|$, where $S \subseteq [d]$. If $|S| = 1$ then the proposition, with $W_+^S$ in place of $W_+^{[d]}$, holds trivially.

Suppose that $2 \leq N \leq d$, and that for every non-empty $Q \subset [d]$ of size $|Q| < N$, the proposition holds with $W_+^Q$ in place of $W_+^{[d]}$. We will show that the proposition then holds for $W_+^S$ whenever $S \subseteq [d]$ and $|S| = N$. Thus let $S \subseteq [d]$ have size $|S| = N$, and let $w(n)$ be a stable sequence in $W_+^S$ with limiting partition $(S_1, \ldots, S_r) \in \mathcal{P}_S$. If $r = 1$ then $w(n)$ necessarily converges in $W_+^S$; since $f$ is continuous on $W_+^S$ (with respect to its natural topology), this implies

$$\lim_{n \to \infty} f(w(n)) = f\left( \lim_{n \to \infty} w(n) \right),$$

which conforms to (38). Next consider the non-trivial case $r \geq 2$. Set $Q = S_1 \cup \ldots \cup S_{r-1}$ and $R = S_r$. Then $w(n)|_Q$ is a stable sequence with limiting partition

$(S_1, \ldots, S_{r-1})$, and by the induction hypothesis

$$\lim_{n \to \infty} f(w(n)|_Q) = \overline{f} \left( \bigoplus_{(S_1,\ldots,S_{r-1})} \lim_{n \to \infty} \left( C_i^{S_1 \cup \ldots \cup S_{j-1}} w(n)_i \right)_{i \in S_j} \right).$$

In particular, $f(w(n)|_Q)$ converges, so by Lemma 12 (with $S$ in place of $[d]$), for every $M$ in the range $1 \leq M \leq 2|Q| - 1$, $f_M(w(n))$ converges, and

(39) $$\lim_{n \to \infty} f_M(w(n)) = \lim_{n \to \infty} f_M(w(n)|_Q).$$

Corollary 15 implies that $w(n)^F|_R$ converges in $W_+^R$, since $R = S_r$ is a cell of the limiting partition of $w(n)$. Therefore $f(w(n)^F|_R)$ converges, and by Lemma 16, for $M$ in the range $2|Q| + 1 \leq M \leq 2|S| - 1$, $f_M(w(n))$ converges, and

$$\lim_{n \to \infty} f_M(w(n)) = \lim_{n \to \infty} f_{M-2|S|} \left( (C_i^Q w(n)_i)_{i \in S_r} \right)$$
(40) $$= f_{M-2|S|} \left( \lim_{n \to \infty} (C_i^Q w(n)_i)_{i \in S_r} \right),$$

the latter equality since $f$ is continuous on $W_+^R$. The behaviour of $f_M(w(n))$ for the remaining index $M = 2|S| - 2|S_r|$ is determined by Lemma 17:

(41) $$\lim_{n \to \infty} f_{2|S|-2|S_r|}(w(n)) = 0.$$

The statements (39), (40) and (41) collectively imply that $\lim f(w(n))$ has the form

$$\lim_{n \to \infty} f(w(n)) = \overline{f} \left( \bigoplus_{(S_1,\ldots,S_{r-1})} \lim_{n \to \infty} \left( C_i^{S_1 \cup \ldots \cup S_{j-1}} w(n)_i \right)_{i \in S_j} \right) \oplus$$

$$f \left( \lim_{n \to \infty} \left( C_i^Q w(n)_i \right)_{i \in S_r} \right)$$

$$= \overline{f} \left( \bigoplus_{(S_1,\ldots,S_r)} \lim_{n \to \infty} \left( C_i^{S_1 \cup \ldots \cup S_{j-1}} w(n)_i \right)_{i \in S_j} \right).$$

Thus the proposition holds for $W_+^S$, completing the induction. ∎

**Corollary 19.** *Let $\widetilde{w}(n)$ be a sequence in $\mathcal{B}_\circ^\Lambda$ that converges in $\mathcal{B}^\Lambda$. Then*

(42) $$\lim_{n \to \infty} \overline{f} \circ \pi(\widetilde{w}(n)) = \overline{f} \circ \pi \left( \lim_{n \to \infty} \widetilde{w}(n) \right).$$

*Proof.* Write $\widetilde{w}(n) = \rho(w(n))$, where $w(n) = \pi(\widetilde{w}(n))$ is by definition a stable sequence in $W_+^{[d]}$, and let $(S_1, \ldots, S_r)$ denote the limiting partition of $w(n)$. Now, $\overline{f} \circ \pi(\widetilde{w}(n)) = \overline{f}(w(n))$, and, by Proposition 18, the left-hand side of (42) is

(43) $$\lim_{n \to \infty} \overline{f}(w(n)) = \overline{f} \left( \bigoplus_{(S_1,\ldots,S_r)} \lim_{n \to \infty} \left( C_i^{S_1 \cup \ldots \cup S_{j-1}} w(n)_i \right)_{i \in S_j} \right).$$

On the other hand, by the calculation (24),

$$\pi\left(\lim_{n\to\infty}\widetilde{w}(n)\right) = \bigoplus_{(S_1,\ldots,S_r)} \lim_{n\to\infty}\left(C_j^{S_1\cup\ldots\cup S_{j-1}}w(n)_i\right)_{i\in S_j},$$

so the right-hand sides of (42) and (43) are equal. ∎

It is a small step from the above corollary to our main result.

**Theorem 20.** *The map $\overline{f}\circ\pi : \mathcal{B}^\Lambda \to \overline{\mathcal{J}}^\Lambda$ is a homeomorphism.*

*Proof.* For brevity write $g = \overline{f}\circ\pi$. Note first that $g$ is a bijection since each of its factors $\overline{f}$ and $\pi$ is. To see that $g$ is continuous, let $w(n) \in \mathcal{B}^\Lambda$ be a convergent sequence. It follows from Corollary 19 and the fact that $\mathcal{B}^\Lambda = \text{cl}\,\mathcal{B}_\circ^\Lambda$ (Proposition 10) that there exists a convergent sequence $w'(n) \in \mathcal{B}_\circ^\Lambda$ such that

(44) $$\lim_{n\to\infty} w'(n) = \lim_{n\to\infty} w(n) \quad \text{and} \quad \lim_{n\to\infty} \|g(w'(n)) - g(w(n))\| = 0.$$

By Corollary 19,

$$\lim_{n\to\infty} g(w'(n)) = g\left(\lim_{n\to\infty} w'(n)\right),$$

which, by (44), implies

$$\lim_{n\to\infty} g(w(n)) = g\left(\lim_{n\to\infty} w(n)\right).$$

Thus $g$ is a continuous bijection; since $\mathcal{B}^\Lambda$ is compact it follows that $g$ is a homeomorphism. ∎

2.8. **The solution to Problem 1.** With Theorem 20 we have a complete solution to Problem 1. Since the composition

$$\mathcal{B}^\Lambda \xrightarrow{\pi} \overline{\mathcal{D}}^\Lambda \xrightarrow{\overline{f}} \overline{\mathcal{J}}^\Lambda$$

is a homeomorphism, the topology that $\overline{\mathcal{J}}^\Lambda$ induces on $\overline{\mathcal{D}}^\Lambda$ via $\overline{f}$ is precisely the topology that $\mathcal{B}^\Lambda$ induces on $\overline{\mathcal{D}}^\Lambda$ via $\pi$. To put it succinctly, the blow-up $(\mathcal{B}^\Lambda, \pi)$ itself constitutes a solution to Problem 1. A sequence $w(n) \in \overline{\mathcal{D}}^\Lambda$ converges in the induced topology if and only if $\pi^{-1}(w(n))$ converges in $\mathcal{B}^\Lambda$, and in this case the limiting value may be computed as

(45) $$\pi\left(\lim_{n\to\infty} \pi^{-1}(w(n))\right),$$

which does not require evaluation of $\overline{f}$. The point is that $\pi$ and $\rho = \pi^{-1}$, unlike $\overline{f}$, are easy to calculate by hand.

More explicitly, a sequence $w(n)$ in the principal component $W_+^{[d]}$ of $\overline{\mathcal{D}}^\Lambda$ is convergent if and only if it is stable. (The characterization of convergent sequences in the other, lower-dimensional pieces of $\overline{\mathcal{D}}^\Lambda$,

$$\bigoplus_{(S_1,\ldots,S_r)} w^{S_i}(n) \in W_+^{(S_1,\ldots,S_r)},$$

is virtually the same: each $w^{S_i}(n)$ must be stable in $W_+^{S_i}$.) Recall that the essence of the computation (45) is captured by the limiting values of moment curves in $W_+^{[d]}$, as in Proposition 9. In distributional notation, given $d\alpha_1, \ldots, d\alpha_r \in \mathcal{D}$ such that

$$(\operatorname{supp} d\alpha_1, \ldots, \operatorname{supp} d\alpha_r) \in \mathcal{P}_\Lambda,$$

the moment curve $d\alpha(t) = d\alpha_1 + t\, d\alpha_2 + \cdots + t^{r-1} d\alpha_r$ tends to the point

$$\begin{aligned} \pi\left(\lim_{t\to 0^+} \pi^{-1}(d\alpha(t))\right) &= \overline{f}^{-1}\left(\lim_{t\to 0^+} f(d\alpha(t))\right) \\ &= d\alpha_1 \oplus p^2_{d\alpha_1} d\alpha_2 \oplus \cdots \oplus p^2_{d\alpha_1+\ldots+d\alpha_{r-1}} d\alpha_r. \end{aligned}$$

The example discussed in Section 2.1 is a simple instance of this general limiting behaviour.

## 3. Isospectral sets of tridiagonal matrices

We turn now to a particular application of the material in Section 2: deducing the topology of the set of symmetric, tridiagonal matrices having prescribed spectrum $\Lambda$. This results in some rather interesting discrete geometry. As mentioned in the Introduction, the isospectral sets we consider have been studied before—our purpose is to provide a new approach.

We state precisely the problem of interest in Section 3.1. Then in Section 3.2 we analyze the combinatorial structure of the semialgebraic set $\mathcal{B}^\Lambda$, constructed earlier in Section 2. We show that $\mathcal{B}^\Lambda$ has a natural "facial" structure, corresponding to the individual pieces of $\overline{\mathcal{D}}^\Lambda$, and that the facial lattice of $\mathcal{B}^\Lambda$ is represented by the partial order $(\mathcal{P}_\Lambda, \preceq)$. The combinatorial structure of $\mathcal{B}^\Lambda$ allows one to very easily construct a polyhedral complex, which we label $\overline{P}_d$, over the manifold of symmetric, tridiagonal matrices that have spectrum $\Lambda$. We describe this construction in Sections 3.3 and 3.4. Finally in Section 3.5 we carry out some explicit computations. We explore the particular example of $\overline{P}_3$: topologically it is a two-holed torus, and its Petrie polygon uses every edge. Concerning the sequence of complexes $\overline{P}_d$, we prove the remarkable fact that hyperbolic tangent generates the Euler characteristics.

We follow the nomenclature for convex polytopes given in Brøndsted's introductory text [2]. In particular we use the characterization of a proper face of a $d$-polytope $P$ in $\mathbf{R}^d$ as the intersection of $P$ with a supporting hyperplane in $\mathbf{R}^d$.

3.1. **Symmetric, tridiagonal matrices; Problem 2.** Let $\mathcal{T}$ denote the set of finite-dimensional, real, symmetric, tridiagonal matrices. We consider the following problem.

**Problem 2** Let $\Lambda = \{\lambda_1, \ldots, \lambda_d\} \subset \mathbf{R}$ be a $d$-element set. Describe the topology of the isospectral set
$$\mathcal{T}^\Lambda = \{T \in \mathcal{T} \mid \operatorname{spec} T = \Lambda\}.$$

The set $\overline{\mathcal{J}}^\Lambda \subset \mathcal{T}^\Lambda$ consists of those matrices whose next-to-diagonal entries are non-negative, and so the parametrization
$$\overline{f} \circ \pi : \mathcal{B}^\Lambda \to \overline{\mathcal{J}}^\Lambda$$
established in Section 2 applies directly to Problem 2. Indeed let $\mathcal{E}$ denote the set of $d$-vectors of the form
$$\varepsilon = (\varepsilon_1, \ldots, \varepsilon_d), \quad \text{where each } \varepsilon_j \in \{+1, -1\}.$$
For each $\varepsilon \in \mathcal{E}$ let $s_\varepsilon$ denote the operation on $d \times d$ matrices of conjugation by $\operatorname{diag}(\varepsilon)$, so that
$$s_\varepsilon(M) = \operatorname{diag}(\varepsilon) M \operatorname{diag}(\varepsilon).$$
Then each map of the form
$$s_\varepsilon \circ \overline{f} \circ \pi : \mathcal{B}^\Lambda \to \mathcal{T}^\Lambda$$
is a homeomorphism from $\mathcal{B}^\Lambda$ onto a "patch" of $\mathcal{T}^\Lambda$, and collectively the maps $s_\varepsilon \circ \overline{f} \circ \pi$ cover $\mathcal{T}^\Lambda$. (Note that the $s_\varepsilon$ are a way of encoding the modification of Definition 1 which would correspond to a particular choice of positive or negative square root for each off-diagonal component of $f$—as Definition 1 stands, only the positive square root is used.) The topology of $\mathcal{T}^\Lambda$ is thus captured by the complex $(\mathcal{B}^\Lambda \times \mathcal{E})/\sim$, consisting of copies $\mathcal{B}^\Lambda \times \{\varepsilon\}$ of $\mathcal{B}^\Lambda$, glued together according to the relation
$$(w, \varepsilon) \sim (w', \varepsilon') \iff s_\varepsilon \circ \overline{f} \circ \pi(w) = s_{\varepsilon'} \circ \overline{f} \circ \pi(w').$$
This gluing map respects the inherent "facial" structure of $\mathcal{B}^\Lambda$, which we describe in Section 3.2 below. In principle, the complex $(\mathcal{B}^\Lambda \times \mathcal{E})/\sim$ is all we need in order to calculate the standard topological invariants of $\mathcal{T}^\Lambda$. However, to make the combinatorial structure of $(\mathcal{B}^\Lambda \times \mathcal{E})/\sim$ more transparent, we go a step further in Section 3.3 and replace $\mathcal{B}^\Lambda$ with a particular convex polytope, call a permutahedron, which is combinatorially equivalent to $\mathcal{B}^\Lambda$ but nevertheless a simpler object. The final result is a polyhedral complex over $\mathcal{T}^\Lambda$ which not only encodes the topology of $\mathcal{T}^\Lambda$, but which is itself rather interesting.

3.2. **The facial structure of $\mathcal{B}^\Lambda$.** For each ordered partition $(S_1, \ldots, S_r) \in \mathcal{P}_{[d]}$, set
$$\mathcal{B}^\Lambda_{(S_1, \ldots, S_r)} = \operatorname{cl}\left(\pi^{-1}\left(W_+^{(S_1, \ldots, S_r)}\right)\right).$$

We refer to $\mathcal{B}^\Lambda_{(S_1,\ldots,S_r)}$ as a face of $\mathcal{B}^\Lambda$, in analogy with the terminology for convex polytopes. We give partial justification for this later in the present section, by showing that, in the appropriate affine setting, each $\mathcal{B}^\Lambda_{(S_1,\ldots,S_r)}$ is the intersection of $\mathcal{B}^\Lambda$ with the bounding hyperplane of a supporting half-space of $\mathcal{B}^\Lambda$, exactly as for a convex polytope. But $\mathcal{B}^\Lambda$ is of course not a polytope; it's faces are in general curved. In the next section we show that the "curved polytope" $\mathcal{B}^\Lambda$ can be flattened, that is, mapped to a bona fide convex polytope that has the same facial structure. Our main objective in this section is to ascertain the inclusion relation between faces of $\mathcal{B}^\Lambda$, and, more precisely, to show that $(\mathcal{P}_{[d]}, \preceq)$ represents the facial lattice of $\mathcal{B}^\Lambda$. This entails what is in essence a technically messier version of the proof of Proposition 10 (which states that $\mathcal{B}^\Lambda_{[d]} = \mathcal{B}^\Lambda$).

To begin, we note the following extension of Proposition 9.

**Proposition 21.** *Let $S \subseteq [d]$, $(S_1,\ldots,S_r) \in \mathcal{P}_S$, and $w^{(S_1,\ldots,S_r)} \in W_+^{(S_1,\ldots,S_r)}$. Then there exists a moment curve $w(t) \in W_+^S$ such that*

$$\overline{f}^{-1}\left(\lim_{t \to 0^+} f(w(t))\right) = w^{(S_1,\ldots,S_r)}.$$

*Proof.* By Theorem 20 the topology that $\overline{\mathcal{J}}^\Lambda$ induces on $\overline{\mathcal{D}}^\Lambda$ via $\overline{f}$ agrees with the topology induced on $\overline{\mathcal{D}}^\Lambda$ by $\mathcal{B}^\Lambda$ via $\pi$. In the case $S = [d]$, Propositions 9 and 10 therefore imply the present proposition directly. But $\Lambda$, and in particular $d = |\Lambda|$, is arbitrary, so once the proposition holds for $S = [d]$, it holds for all $S \subset [d]$, since we may replace $\Lambda$ with $\Lambda' = \{\lambda_n \in \Lambda \mid n \in S\}$. ∎

For each $(S_1,\ldots,S_r) \in \mathcal{P}_{[d]}$, let $F_{(S_1,\ldots,S_r)} \subset \mathcal{B}^\Lambda$ denote the set

$$F_{(S_1,\ldots,S_r)} = \pi^{-1}\left(\bigcup_{(R_1,\ldots,R_q) \preceq (S_1,\ldots,S_r)} W_+^{(R_1,\ldots,R_q)}\right).$$

**Lemma 22.** *For each $(S_1,\ldots,S_r) \in \mathcal{P}_{[d]}$, $F_{(S_1,\ldots,S_r)} \subseteq \mathcal{B}^\Lambda_{(S_1,\ldots,S_r)}$.*

*Proof.* Let $\widetilde{w} \in F_{(S_1,\ldots,S_r)}$, and write $w^{(R_1,\ldots,R_q)} = \pi(\widetilde{w})$, so that $(R_1,\ldots,R_q) \preceq (S_1,\ldots,S_r)$ and hence has the form

$$(R_1,\ldots,R_q) = (R_1^1,\ldots,R_{i_1}^1,\ldots,R_1^r,\ldots,R_{i_r}^r),$$

where each $(R_1^j,\ldots,R_{i_j}^j) \in \mathcal{P}_{S_j}$. Using Proposition 21, for each $S_j$, let $w^{S_j}(t) \in W_+^{S_j}$ be a moment curve such that

$$\overline{f}^{-1}\left(\lim_{t \to 0^+} f(w^{S_j}(t))\right) = w^{(R_1^j,\ldots,R_{i_j}^j)},$$

where $w^{(R_1^j,\ldots,R_{i_j}^j)}$ denotes the appropriate component of $w^{(R_1,\ldots,R_q)}$. Then $\bigoplus_{(S_1,\ldots,S_r)} w^{S_j}(t)$ is a moment curve in $W_+^{(S_1,\ldots,S_r)}$, and

$$\begin{aligned}
\lim_{t\to 0^+} \overline{f}\left(\bigoplus_{(S_1,\ldots,S_r)} w^{S_j}(t)\right) &= \lim_{t\to 0^+} \bigoplus_{(S_1,\ldots,S_r)} f(w^{S_j}(t)) \\
&= \bigoplus_{(S_1,\ldots,S_r)} \lim_{t\to 0^+} f(w^{S_j}(t)) \\
&= \bigoplus_{(S_1,\ldots,S_r)} \overline{f}(w^{(R_1^j,\ldots,R_{i_j}^j)}) \\
&= \overline{f}(w^{(R_1,\ldots,R_q)}),
\end{aligned}$$

whereby

$$\overline{f}^{-1}\left(\lim_{t\to 0^+} \overline{f}\left(\bigoplus_{(S_1,\ldots,S_r)} w^{S_j}(t)\right)\right) = w^{(R_1,\ldots,R_q)}.$$

By Theorem 20, this implies

$$\lim_{t\to 0^+} \pi^{-1}\left(\bigoplus_{(S_1,\ldots,S_r)} w^{S_j}(t)\right) = \widetilde{w},$$

so that $\widetilde{w} \in \mathcal{B}_{(S_1,\ldots,S_r)}^\Lambda = \mathrm{cl}\left(\pi^{-1}\left(W_+^{(S_1,\ldots,S_r)}\right)\right)$. ∎

In Section 2 it was convenient to have $\mathcal{B}^\Lambda$ in projective space so that we could use homogeneous coordinates, and not have to worry about rescaling. For present considerations, it is better to embed $\mathcal{B}^\Lambda$ in affine space. To that end, for each non-empty $S \subseteq [d]$, let

$$\mathbf{R}^S = \{w^S : S \to \mathbf{R}\}$$

be the $|S|$-dimensional Euclidean space of real-valued sequences indexed over $S$, and set

$$E = \prod_{\emptyset \neq S \subseteq [d]} \mathbf{R}^S.$$

Let $A$ denote the affine subspace of $E$ defined as

$$A = \left\{w \in E \;\middle|\; \forall S \subseteq [d], S \neq \emptyset, \sum_{n \in S} w_n^S = 1\right\}.$$

Now, the ambient space of $\mathcal{B}^\Lambda$, $\mathcal{W} = \prod_{\emptyset \neq S \subseteq [d]} W^S$, is a quotient space of $E - \{0\}$, and the inequalities (17) defining $\mathcal{B}^\Lambda$ guarantee that each $w \in \mathcal{B}^\Lambda$ has a unique normalized representor $\widehat{w} \in A$. Hence we may regard $\mathcal{B}^\Lambda \subset A$ as sitting in affine space, without altering its topological structure; we do so for the remainder of the present section.

We now define a family of half-spaces $H^+_{(S_1,\ldots,S_r)}$ of $A$, along with their bounding hyperplanes $H_{(S_1,\ldots,S_r)}$. It follows from the the explicit construction of $\rho = \pi^{-1}$ given in Definition 6 that for $w \in \pi^{-1}\left(W^{(S_1,\ldots,S_r)}_+\right)$, the set
$$\mathcal{I} = \left\{(n, S) \,\big|\, n \in S \subseteq [d],\, w^S_n = 0\right\}$$
does not depend on the particular point $w$. Rather,
$$\mathcal{I} = \mathcal{I}_{(S_1,\ldots,S_r)}$$
is determined by the partition $(S_1,\ldots,S_r) \in \mathcal{P}_{[d]}$. Moreover,
$$\mathcal{I}_{(S_1,\ldots,S_r)} \subseteq \mathcal{I}_{(R_1,\ldots,R_q)} \iff (R_1,\ldots,R_q) \preceq (S_1,\ldots,S_r). \tag{46}$$
Now, given $(S_1,\ldots,S_r) \in \mathcal{P}_{[d]}$, let
$$H_{(S_1,\ldots,S_r)} \subset A$$
denote the hyperplane in $A$ determined by the equation
$$\sum_{(n,S) \in \mathcal{I}_{(S_1,\ldots,S_r)}} w^S_n = 0,$$
and let $H^+_{(S_1,\ldots,S_r)} \subset A$ denote the corresponding half-space, determined by the inequality
$$\sum_{(n,S) \in \mathcal{I}_{(S_1,\ldots,S_r)}} w^S_n \geq 0.$$

**Lemma 23.** *For every $(S_1,\ldots,S_r) \in \mathcal{P}_{[d]}$, the half-space $H^+_{(S_1,\ldots,S_r)}$ includes $\mathcal{B}^\Lambda$, and $F_{(S_1,\ldots,S_r)} = H_{(S_1,\ldots,S_r)} \cap \mathcal{B}^\Lambda$.*

*Proof.* The inequalities (17) defining $\mathcal{B}^\Lambda$ and the defining equations of $A$ imply that if $w \in \mathcal{B}^\Lambda$, then its normalized representor $\widehat{w} \in A$ satisfies $\widehat{w}^S_n \geq 0$ for every $n \in S \subseteq [d]$. It follows that $\mathcal{B}^\Lambda \subset H^+_{(S_1,\ldots,S_r)}$.

Also, $\widehat{w} \in H_{(S_1,\ldots,S_r)}$ if and only if for each $(n,S) \in \mathcal{I}_{(S_1,\ldots,S_r)}$, $\widehat{w}^S_n = 0$. By (46) this condition is equivalent to the condition that
$$\pi(\widehat{w}) \in W^{(R_1,\ldots,R_q)}_+,$$
for some $(R_1,\ldots,R_q) \preceq (S_1,\ldots,S_r)$, or in other words, $\widehat{w} \in F_{(S_1,\ldots,S_r)}$. ∎

Combining Lemmas 22 and 23 yields the following characterization of the faces $\mathcal{B}^\Lambda_{(S_1,\ldots,S_r)}$.

**Proposition 24.** *For every $(S_1,\ldots,S_r) \in \mathcal{P}_{[d]}$, $\mathcal{B}^\Lambda_{(S_1,\ldots,S_r)} = F_{(S_1,\ldots,S_r)}$.*

*Proof.* By definition, $F_{(S_1,\ldots,S_r)} \supseteq \pi^{-1}\left(W^{(S_1,\ldots,S_r)}_+\right)$. Using Lemma 23, $F_{(S_1,\ldots,S_r)} = \mathcal{B}^\Lambda \cap H_{(S_1,\ldots,S_r)}$ is closed, so
$$\mathcal{B}^\Lambda_{(S_1,\ldots,S_r)} = \mathrm{cl}\left(\pi^{-1}\left(W^{(S_1,\ldots,S_r)}_+\right)\right) \subseteq F_{(S_1,\ldots,S_r)}.$$
Lemma 22 gives the reverse inclusion. ∎

**Corollary 25.** $\mathcal{B}^\Lambda_{(R_1,\ldots,R_q)} \subseteq \mathcal{B}^\Lambda_{(S_1,\ldots,S_r)}$ *if and only if* $(R_1,\ldots,R_q) \preceq (S_1,\ldots,S_r)$.

Thus we have shown that $(\mathcal{P}_{[d]}, \preceq)$ represents the facial lattice of $\mathcal{B}^\Lambda$. Note that the dimension of the face $\mathcal{B}^\Lambda_{(S_1,\ldots,S_r)}$ is $\dim(W_+^{(S_1,\ldots,S_r)}) = d - r$. The least-dimensional faces, corresponding to ordered partitions of the form

$$(\{\sigma(1)\},\ldots,\{\sigma(d)\}),$$

where $\sigma$ is a permutation of $[d]$, consist of a single point and have dimension 0. If two ordered partitions

$$(R_1,\ldots,R_q) \preceq (S_1,\ldots,S_r)$$

are adjacent with respect to $\preceq$, that is if $r = q+1$, then evidently

$$\dim\left(\mathcal{B}^\Lambda_{(S_1,\ldots,S_r)}\right) = \dim\left(\mathcal{B}^\Lambda_{(R_1,\ldots,R_q)}\right) + 1,$$

a fact which we will use in the next section.

3.3. **Permutahedra.** Permutahedra are convex polytopes that have a very special structure. For instance, they are simple zonotopes that tile Euclidean space; see Ziegler [11, pp.17,18,200]. We give a slightly different description from that of Ziegler, and consequently our characterization of faces, Proposition 27, differs somewhat from his. But it is easy to reconcile the two descriptions, so we do not bother to prove Proposition 27.

We fix notation as follows. Let $S_{[d]}$ denote the symmetric group, consisting of all permutations of $[d] = \{1,\ldots,d\}$. For each non-empty subset $R \subseteq [d]$, let $S_R$ denote the subgroup of $S_{[d]}$ consisting of those permutations which fix each element of $[d]-R$; let $\mathcal{G}$ denote the collection of subgroups of $S_{[d]}$ given by

$$\mathcal{G} = \left\{ S_{R_1} \times \cdots \times S_{R_q} \mid (\{1\},\ldots,\{d\}) \preceq (R_1,\ldots,R_q) \in \mathcal{P}_{[d]} \right\}.$$

We say that a non-empty set $R \subseteq [d]$ is an invariant set of a group $H \subseteq S_{[d]}$ if every $\sigma \in H$ fixes $R$ setwise; if in addition no non-empty proper subset of $R$ is invariant, we say $R$ is a minimal invariant set of $H$. Note that for

$$G = S_{R_1} \times \cdots \times S_{R_q} \in \mathcal{G},$$

the sets $R_1,\ldots,R_q$ are the minimal invariant subsets of $G$, and they have a unique ordering for which $(\{1\},\ldots,\{d\}) \preceq (R_1,\ldots,R_q) \in \mathcal{P}_{[d]}$. The collection $\mathcal{G}$ will be useful in discussing the structure of the permutahedron, which we now define.

**Definition 26.** *Let $P_d$ denote the convex hull of the set*

(47) $$\left\{ (\sigma(1),\ldots,\sigma(d)) \in \mathbf{R}^d \mid \sigma \in S_{[d]} \right\}.$$

*The convex polytope $P_d$ is called a permutahedron.*

Given a set $F \subset \mathbf{R}^d$, let $\widehat{F}$ denote its convex hull. We will generally identify elements of $S_{[d]}$ with points in $\mathbf{R}^d$ as in (47); with this convention we have, for instance,

$$P_d = \widehat{S_{[d]}}.$$

The set of groups $\mathcal{G}$ provides a simple description of the faces of $P_d$, as follows. Let $\mathcal{G}S_{[d]}$ denote the set of right cosets of the form $G\sigma$, where $G \in \mathcal{G}$ and $\sigma \in S_{[d]}$.

**Proposition 27.** *The faces of $P_d$ are precisely the sets $\widehat{F}$ where $F \in \mathcal{G}S_{[d]}$.*

The above proposition implies that
$$\left(\mathcal{G}S_{[d]}, \subseteq\right)$$
represents the facial lattice of $P_d$. Indeed this corresponds to the facial lattice of $\mathcal{B}^\Lambda$ via the following map.

**Definition 28.** *We denote by $\varphi$ the map*
$$\varphi : \mathcal{G}S_{[d]} \to \mathcal{P}_{[d]},$$
*constructed as follows. Let $H \in \mathcal{G}S_{[d]}$. Choose $\sigma \in H$ arbitrarily; then $H\sigma^{-1} \in \mathcal{G}$. Let $R_1, \ldots, R_q \subset [d]$ be the minimal invariant sets of $H\sigma^{-1}$, ordered such that*
$$(\{1\}, \ldots, \{d\}) \preceq (R_1, \ldots, R_q).$$
*Set*
$$\varphi(H) = \left(\sigma^{-1}(R_1), \ldots, \sigma^{-1}(R_q)\right) \in \mathcal{P}_{[d]}.$$

Note that if $H \in \mathcal{G}S_{[d]}$ corresponds to a vertex of $P_d$, that is, if $H = \{\sigma\}$ is a singleton, then according to the above definition
$$\varphi(H) = \left(\{\sigma^{-1}(1)\}, \ldots, \{\sigma^{-1}(d)\}\right).$$
In other words, on the level of vertices, $\varphi$ coincides essentially with the map $\sigma \mapsto \sigma^{-1}$. More generally, let $H \in \mathcal{G}S_{[d]}$ and write $(T_1, \ldots, T_r) = \varphi(H)$. Then it is a straightforward to verify that
$$(48) \qquad H = \left\{\sigma \in S_{[d]} \mid (\{1\}, \ldots, \{d\}) \preceq (\sigma(T_1), \ldots, \sigma(T_r))\right\},$$
whence $\varphi$ is one-to-one. And given any $(T_1, \ldots, T_r) \in \mathcal{P}_{[d]}$, the set $H$ defined by (48) belongs to $\mathcal{G}S_{[d]}$, and so $\varphi$ is onto; that is, $\varphi$ is a bijection. Indeed, $P_d$ and $\mathcal{B}^\Lambda$ are combinatorially equivalent, in the following sense.

**Proposition 29.** *The map $\varphi : \mathcal{G}S_{[d]} \to \mathcal{P}_{[d]}$ is a lattice isomorphism, so that*
$$\left(\mathcal{G}S_{[d]}, \subseteq\right) \cong \left(\mathcal{P}_{[d]}, \preceq\right).$$

*Proof.* Since $\varphi$ is a bijection, we need only show that for every $H, K \in \mathcal{G}S_{[d]}$, $H \subseteq K$ if and only if $\varphi(H) \preceq \varphi(K)$. To that end, let $H, K \in \mathcal{G}S_{[d]}$ and write $(R_1, \ldots, R_q) = \varphi(H)$, $(T_1, \ldots, T_r) = \varphi(K)$. By (48) the condition $H \subseteq K$ is evidently equivalent to the condition that for every $\sigma \in S_{[d]}$,
$$(\{1\}, \ldots, \{d\}) \preceq (\sigma(R_1), \ldots, \sigma(R_q)) \Rightarrow (\{1\}, \ldots, \{d\}) \preceq (\sigma(T_1), \ldots, \sigma(T_r)).$$
And the latter condition is easily shown to be equivalent to the relation
$$(R_1, \ldots, R_q) \preceq (T_1, \ldots, T_r).$$
∎

For an arbitrary convex polytope $P$, if two faces $F \subset F'$ are adjacent in the facial lattice of $P$, that is, if there is no face $F''$ of $P$ such that $F \subset F'' \subset F'$, then $\dim(F)' = \dim(F) + 1$. As pointed out in the previous section, the facial lattice of $\mathcal{B}^\Lambda$ also has this property; therefore the map $\varphi$ automatically respects dimension:

$$\dim(\widehat{F}) = \dim(\mathcal{B}^\Lambda_{\varphi(F)}),$$

for every $F \in \mathcal{GS}_{[d]}$.

Any isomorphism $\iota$ between the facial lattices of two convex polytopes $P$ and $P'$ has a topological realization, in the sense that there exists a homeomorphism $h : P \to P'$ which preserves faces, and whose action on faces agrees with $\iota$. By the same token, the map $\varphi$ has a topological realization

$$h : P_d \to \mathcal{B}^\Lambda.$$

We sketch briefly one possible construction of such an $h$. Here it is useful to have a metric on $\mathcal{B}^\Lambda$, so that we may speak of geodesics. As in the previous section, let us regard $\mathcal{B}^\Lambda$ as sitting in the affine space $A$, and consider the metric that $\mathcal{B}^\Lambda$ inherits from the ambient Euclidean space $E = \prod_{\emptyset \neq S \subseteq [d]} \mathbf{R}^{|S|}$.

For each face $\mathcal{B}^\Lambda_{(S_1,\ldots,S_r)}$ of $\mathcal{B}^\Lambda$, consider the point $w^{(S_1,\ldots,S_r)} \in W_+^{(S_1,\ldots,S_r)}$ having constant entries $w_n^{S_j} = 1$ for every $j$, $n \in S_j$. We refer to the normalized representor $b_{(S_1,\ldots,S_r)} \in A$ of $\pi^{-1}(w^{(S_1,\ldots,S_r)})$ as the "barycentre" of $\mathcal{B}^\Lambda_{(S_1,\ldots,S_r)}$. The map $h$ may be constructed recursively as follows.

<u>Initial step</u>  For each vertex $\sigma \in P_d$, set $h(\sigma) = \mathcal{B}^\Lambda_{\varphi(\{\sigma\})}$ (the latter consists of just a single point).

<u>Continuing step</u>  Suppose that $h$ has been defined on the $j$-skeleton of $P_d$, that is, on the set of faces of $P_d$ having dimension at most $j$, where $0 \leq j \leq d-2$. Given a $(j+1)$-dimensional face $\widehat{F}$, extend $h$ to $\widehat{F}$ as follows. Let $h$ map the barycentre $f_0$ of $\widehat{F}$ (the mean of $F$) to the barycentre $b_{\varphi(F)}$ of $\mathcal{B}^\Lambda_{\varphi(F)}$. And for each point $f_1$ in the $j$-skeleton of $\widehat{F}$, let $h$ map the segment $\overline{f_0 f_1}$ to the geodesic in $\mathcal{B}^\Lambda_{(S_1,\ldots,S_r)}$ connecting $b_{\varphi(F)}$ to $h(f_1)$, at constant speed with respect to arclength.

The map $h : P_d \to \mathcal{B}^\Lambda$ so constructed is a homeomorphism, and its inverse can be viewed as a "flattening" of $\mathcal{B}^\Lambda$. The permutahedron $P_d$ is homeomorphic to $\overline{\mathcal{J}}^\Lambda$ via the composite map

$$(49) \qquad P_d \xrightarrow{h} \mathcal{B}^\Lambda \xrightarrow{\pi} \overline{\mathcal{D}}^\Lambda \xrightarrow{\overline{f}} \overline{\mathcal{J}}^\Lambda.$$

This serves as the basis for a polyhedral complex, which we describe in the next section, on the set of $d \times d$ symmetric, tridiagonal matrices having spectrum $\Lambda$.

3.4. **A polyhedral complex over $\mathcal{T}^\Lambda$.** Referring to the composite map (49), for $p \in P_d$ write

$$J_p = \overline{f} \circ \pi \circ h(p).$$

Let $\mathcal{E} \subset \mathbf{R}^d$ denote the $2^d$-element set consisting of maps $\varepsilon$ of the form

$$\varepsilon : [d] \to \{+1, -1\}.$$

(In other words, $\mathcal{E}$ consists of the vertices of the twice-unit cube centred at the origin.) For each $(p, \varepsilon) \in P_d \times \mathcal{E}$, set

$$g(p, \varepsilon) = \operatorname{diag}(\varepsilon) J_p \operatorname{diag}(\varepsilon).$$

Since the map $p \mapsto J_p$ is a homeomorphism from $P_d$ onto $\overline{\mathcal{J}}^\Lambda$, the map

$$g : P_d \times \mathcal{E} \to \mathcal{T}^\Lambda$$

is easily seen to be a continuous surjection. Let $\sim$ denote the relation defined on $P_d \times \mathcal{E}$ defined by

$$(p, \varepsilon) \sim (p', \varepsilon') \iff g(p, \varepsilon) = g(p', \varepsilon').$$

And let $[(p, \varepsilon)]$ denote the equivalence class of $(p, \varepsilon)$ modulo $\sim$. The relation $\sim$ may be characterized purely in terms of $P_d \times \mathcal{E}$ as follows.

**Proposition 30.** *Let $(p, \varepsilon), (p', \varepsilon') \in P_d \times \mathcal{E}$. Let $F \in \mathcal{G}S_{[d]}$ represent the unique face $\widehat{F}$ of $P_d$ whose relative interior contains $p$. Choose $\sigma \in F$, and write $R_1, \ldots, R_q$ for the minimal invariant sets of the group $F\sigma^{-1} \in \mathcal{G}$, ordered such that*

$$(\{1\}, \ldots, \{d\}) \preceq (R_1, \ldots, R_q).$$

*Then $(p, \varepsilon) \sim (p', \varepsilon')$ if and only if $p = p'$ and $\varepsilon|_{R_j} = \pm \varepsilon'|_{R_j}$ for each $R_j$.*

*Proof.* The matrix $J_p$ may be obtained from $g(p, \varepsilon)$ by replacing off-diagonal entries with their absolute values. Since $p \mapsto J_p$ is one-to-one, it follows that $g(p, \varepsilon)$ determines $p$ uniquely. Thus $(p, \varepsilon) \sim (p', \varepsilon')$ implies $p = p'$ and we need only characterize the relation $(p, \varepsilon) \sim (p, \varepsilon')$.

The fact that $p$ is in the relative interior of $\widehat{F}$ implies that $h(p)$ is in the relative interior of $\mathcal{B}^\Lambda_{\varphi(F)}$, or in other words, in $\pi^{-1}\left(W^{\varphi(F)}_+\right)$. Therefore $J_p \in \overline{f}\left(W^{\varphi(F)}_+\right)$, which, recalling that $\varphi(F) = (\sigma^{-1}(R_1), \ldots, \sigma^{-1}(R_q))$, implies

$$J_p = J_1 \oplus \cdots \oplus J_q,$$

where each $J_i$ is an $|R_i| \times |R_i|$ Jacobi matrix. With this decomposition, the condition $(p, \varepsilon) \sim (p, \varepsilon')$ becomes

$$\operatorname{diag}(\varepsilon) J_1 \oplus \cdots \oplus J_q \operatorname{diag}(\varepsilon) = \operatorname{diag}(\varepsilon') J_1 \oplus \cdots \oplus J_q \operatorname{diag}(\varepsilon'),$$

which is evidently equivalent to the condition that, for each $R_i$, $\varepsilon|_{R_i} = \pm \varepsilon'|_{R_i}$. ∎

We now define the polyhedral complex of interest.

**Definition 31.** *Let $\overline{P}_d$ denote the quotient $(P_d \times \mathcal{E})/\sim$, and define*

$$\widetilde{g} : \overline{P}_d \to \mathcal{T}^\Lambda$$

*to be the map induced by $g$, so that $\widetilde{g}([(p, \varepsilon)]) = g(p, \varepsilon)$.*

Faces of $\overline{P}_d$ are sets of the form $\widehat{F} \times \{\varepsilon\}$ where $\varepsilon \in \mathcal{E}$ and $\widehat{F}$ is a face of $P_d$. Proposition 30 implies that either the relative interiors of two faces $\widehat{F} \times \{\varepsilon\}$ and $\widehat{F'} \times \{\varepsilon'\}$ are disjoint, or the two faces are equal. And the latter can only happen if $F = F'$. In the case of a least-dimensional face, a vertex $\{\sigma\}$ where $\sigma \in S_{[d]}$, we have $\{\sigma\} \times \{\varepsilon\} = \{\sigma\} \times \{\varepsilon'\}$ in $\overline{P}_d$ for every $\varepsilon \in \mathcal{E}$. In other words, $\overline{P}_d$ has the same vertices as $P_d$, one for each permutation $\sigma \in S_{[d]}$. At the other extreme, there are $2^{d-1}$ distinct faces of $\overline{P}_d$ of the form $P_d \times \{\varepsilon\}$; a given set $P_d \times \{\varepsilon\}$ gets identified in $\overline{P}_d$ with the set $P_d \times \{\varepsilon'\}$ only when $\varepsilon' = \pm\varepsilon$.

The following result encapsulates the solution to Problem 2. To prove it, we need only verify that the map $\widetilde{g}$ is a homeomorphism.

**Theorem 32.** $\overline{P}_d$ *is a polyhedral complex over* $\mathcal{T}^\Lambda$.

*Proof.* Since $g : P_d \times \mathcal{E} \to \mathcal{T}^\Lambda$ is a continuous surjection, $\widetilde{g}$ is a continuous bijection by construction. That $\widetilde{g}$ is a homeomorphism then follows from compactness of $\overline{P}_d$. ∎

Thus the topology of $\mathcal{T}^\Lambda$ is encoded by the polyhedral complex $\overline{P}_d$. Using Proposition 30, the calculation of topological invariants of $\overline{P}_d$ is reasonably straightforward; we illustrate this in the next section.

### 3.5. Illustration of the results, Euler characteristics.

In this section we apply the earlier results specifically to $3 \times 3$ matrices, and we calculate the Euler characteristic of $\mathcal{T}^\Lambda$ for an arbitrary $d$-element set $\Lambda = \{\lambda_1, \dots, \lambda_d\} \subset \mathbf{R}$.

To begin, consider the permutahedron $P_d$, the basic unit of the complex $\overline{P}_d$, and let us count the number of faces having a fixed dimension $n$, where $0 \leq n \leq d - 1$. A face $\widehat{F}$ of $P_d$ has dimension $n$ if and only if the ordered partition

$$(S_1, \dots, S_r) = \varphi(F)$$

has $r = d - n$ cells. There are precisely

$$(50) \qquad (d-n)!S(d, d-n)$$

such ordered partitions, where $S(d, d-n)$ denotes the Stirling number of the second kind.

Recall that if $F, F' \in \mathcal{G}S_{[d]}$ are distinct, then for any $\varepsilon, \varepsilon' \in \mathcal{E}$ so are the faces $\widehat{F} \times \{\varepsilon\}, \widehat{F'} \times \{\varepsilon\}$ of $\overline{P}_d$. On the other hand, for fixed $F$, Proposition 30 tells us precisely how many faces of $\overline{P}_d$ there are of the form $\widehat{F} \times \{\varepsilon\}$ ($\varepsilon \in \mathcal{E}$). Namely, let $N$ denote the number in question, and write

$$(\sigma^{-1}(R_1), \dots, \sigma^{-1}(R_q)) = \varphi(F),$$

where $\sigma \in F$ and $R_1, \dots, R_q$ are the minimal invariant subsets of $F\sigma^{-1}$. Then $N$ is the number of equivalence classes of $\mathcal{E}$ modulo the relation

$$\varepsilon|_{R_j} = \pm\varepsilon'|_{R_j} \quad \text{for every } R_j.$$

That is,

$$(51) \qquad N = 2^d/2^q = 2^n,$$

where $n = d - q$ is the dimension of $\widehat{F} \times \{\varepsilon\}$. Combining (50) and (51) yields the number of $n$-dimensional faces of $\overline{P}_d$ and hence the following explicit form for the Euler characteristic.

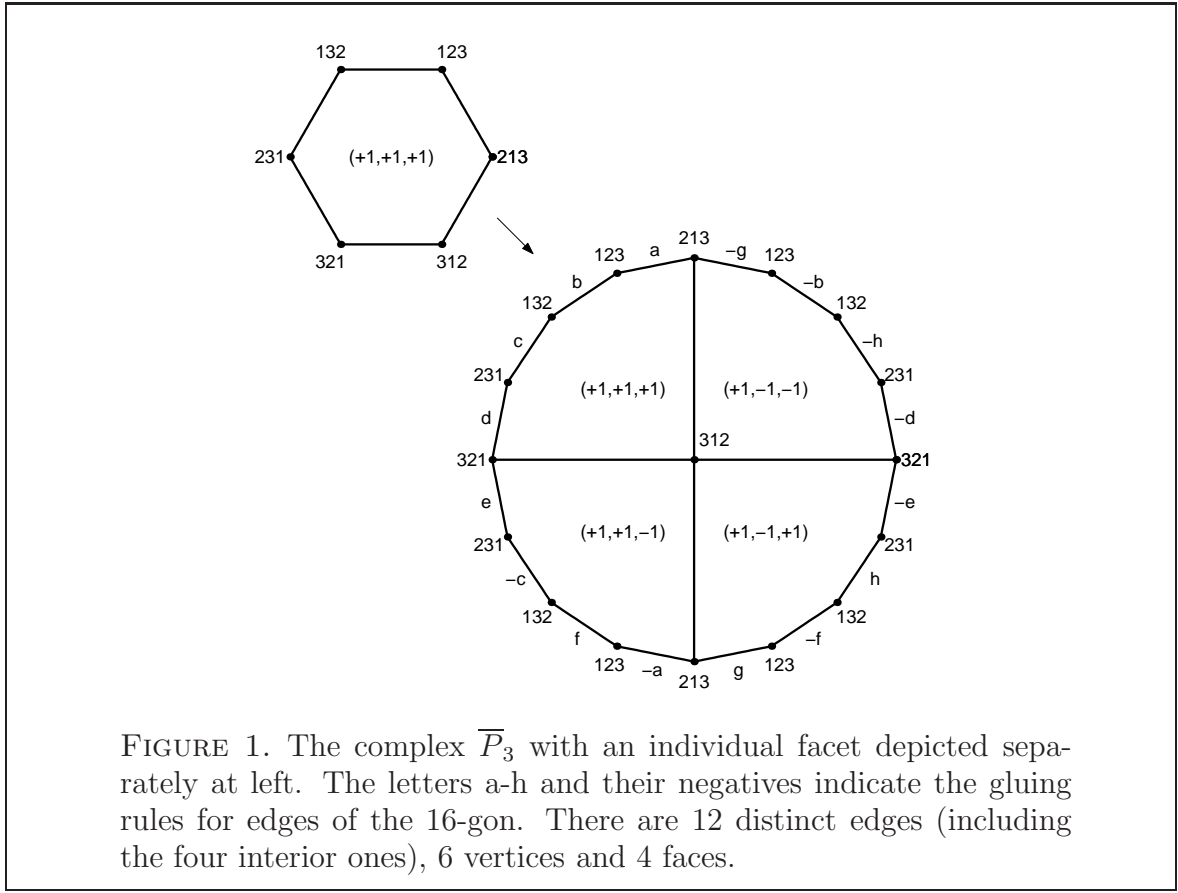

FIGURE 1. The complex $\overline{P}_3$ with an individual facet depicted separately at left. The letters a-h and their negatives indicate the gluing rules for edges of the 16-gon. There are 12 distinct edges (including the four interior ones), 6 vertices and 4 faces.

**Theorem 33.** *The Euler characteristic of $\overline{P}_d$ is*

$$\chi(\overline{P}_d) = \sum_{n=0}^{d-1}(-2)^n(d-n)!S(d,d-n) = \sum_{q=1}^{d}(-2)^{d-q}q!S(d,q).$$

We show below that the above value of the Euler characteristic has a simple generating function. But first we discuss the particular complex $\overline{P}_3$. It consists of $2^{3-1} = 4$ copies of the regular hexagon $P_3$, glued together according to the relation $\sim$ as depicted in Figure 1. By Theorem 33 the Euler characteristic is

$$\chi\left(\overline{P}_3\right) = -2,$$

which, using Theorem 32 and the standard classification of 2-manifolds, implies:

**Theorem 34.** *$\mathcal{T}^\Lambda$ is a two-holed torus whenever $\Lambda$ is a set of 3 distinct real numbers.*

From the point of view of discrete geometry, the polyhedral complexes $\overline{P}_d$ seem to be of intrinsic interest, quite apart from their relationship to tridiagonal matrices. For instance, J. Bokowski has made the following observation. See [4, pp.24,25,223] for relevant definitions and comparison with the Platonic solids.

**Proposition 35.** *The Petrie polygon of $\overline{P}_3$ covers every edge.*

This unusual property can be verified directly in Figure 1. Another remarkable property of the complexes $\overline{P}_d$—and also of the isospectral sets $\mathcal{T}^\Lambda$—is that hyperbolic tangent generates the Euler characteristics.

**Theorem 36.**
$$\tanh x = \sum_{d=1}^{\infty} \chi\left(\overline{P}_d\right) \frac{x^d}{d!}.$$

*Proof.* According to a theorem of Frobenius [3, p. 244], the so-called Eulerian polynomial of degree $d$ has the representation
$$p_d(u) = u \sum_{r=1}^{d} r! S(d,r)(u-1)^{d-r}.$$

It is known [3, p. 259] that
$$\tan x = \sum_{d=1}^{\infty} |p_d(-1)| \frac{x^d}{d!}$$
$$= \sum_{d=1}^{\infty} i^{d+1} p_d(-1) \frac{x^d}{d!},$$
where $i^2 = -1$. It follows from the relation
$$\tanh x = -i \tan(ix)$$
that
$$\tanh x = \sum_{d=1}^{\infty} (-1)^d p_d(-1) \frac{x^d}{d!}$$
$$= \sum_{d=1}^{\infty} -p_d(-1) \frac{x^d}{d!},$$
the latter equality because tanh is odd. Comparing Theorem 33 to the Frobenius representation of $p_d$ shows $-p_d(-1) = \chi(\overline{P}_d)$. ∎

Peter C. Gibson, Department of Mathematics, University of Washington, Seattle, WA 98195-4350, USA